\documentclass[12pt, reqno]{amsart}
\usepackage[utf8]{inputenc}
\usepackage{amsaddr}

\usepackage{import}
\usepackage{graphicx}
\usepackage{color}
\usepackage{transparent}
\usepackage{calc}
\usepackage{float}

\usepackage{subcaption}
\usepackage{overpic}
\usepackage{amsmath,amsthm,amssymb,amscd,amsfonts,amsbsy}
\usepackage{color}
\usepackage{float}
\usepackage{hyperref}
\usepackage{enumerate}
\usepackage{enumitem}    
\usepackage{mathrsfs}
\usepackage{mathtools}
\usepackage{thmtools}
\usepackage{comment}
\usepackage{thm-restate}
\usepackage{cleveref}
\usepackage{graphics, setspace}
\usepackage{listings,xcolor}
\usepackage[normalem]{ulem}
\usepackage[toc,page]{appendix}
\usepackage[top=3cm,bottom=2cm,left=1.5cm,right=1.5cm]{geometry}
\usepackage{soul}
\usepackage{dsfont}   

\usepackage[toc,page]{appendix}

\allowdisplaybreaks

\newcommand{\mathsym}[1]{{}}
\newcommand{\unicode}[1]{{}}

\newcommand{\N}{\ensuremath{\mathbb{N}}}

\newcommand{\R}{\ensuremath{\mathbb{R}}}

\newcommand{\ov}{\overline}

\newcommand{\flow}{\varphi}

\newcommand{\sgn}{\mathrm{sign}}
\newcommand{\set}[1]{\left\{#1\right\}}
\newcommand{\hausdorff}[1]{\operatorname{dim_H}\left(#1\right)}
\newcommand{\cal}[1]{\mathcal{#1}}
\newcommand{\abs}[1]{\left\lvert{#1}\right\rvert}
\newcommand{\Bigabs}[1]{\Big|{#1}\Big|}
\newcommand{\norm}[1]{\left\lVert{#1}\right\rVert}
\newcommand{\inner}[1]{\langle{#1}\rangle}
\newcommand{\cl}[1]{\ensuremath{\overline{#1}}}
\newcommand{\Dom}[1]{\ensuremath{\operatorname{Dom}(#1)}}
\newcommand{\proj}[0]{\ensuremath{\operatorname{proj}}}
\newcommand{\opname}[1]{\ensuremath{\operatorname{#1}}}

\theoremstyle{plain}
\newtheorem{theorem}{Theorem}[section]
\newtheorem{proposition}{Proposition}[section]

\newtheorem{lemma}{Lemma}[section]

\newtheorem{maintheorem}{Theorem}

\theoremstyle{definition}
\newtheorem{definition}{Definition}[section]
\newtheorem{example}{Example}[section]
\newtheorem{remark}{Remark}[section]

\def\R{\mathbb R}
\def\N{\mathbb N}

\def \y {\mathbf{y}}

\renewcommand{\subset}{\subseteq}
\renewcommand{\supset}{\supseteq}

\graphicspath{{figures/}}

\title[On Bykov Cycles in Filippov Systems]{On Bykov Cycles in Filippov Systems}

\author[M. G. C. Cunha, and P. R. Silva]
{Matheus G. C. Cunha$^{1}$, and Paulo R. Silva$^2$}

\address{$^{1,2}$Departamento de Matem\'{a}tica - Instituto de Bioci\^{e}ncias, Letras e Ci\^{e}ncias Exatas (IBILCE) - Universidade
Estadual Paulista (UNESP), \ Rua Crist\'{o}vao Colombo, 2265, Jardim Nazareth, CEP 15054-000, Sao Jos\'{e} do Rio Preto, SP,
Brazil}
\email{$^1$matheus.gc.cunha@unesp.br}
\email{$^2$paulo.r.silva@unesp.br}

\begin{document}
\subjclass[2020]{Primary: 34A36, 34C37; Secondary: 37B10, 37D35, 37G10}

\keywords{Bykov Cycles, Filippov Systems, Sliding Dynamics, Piecewise Expanding Maps, Regularization}

\maketitle

\begin{abstract}
We introduce two types of Bykov cycles for Filippov systems: crossing Bykov cycles, which retain the main dynamical features of their smooth counterparts, and sliding Bykov cycles, whose dynamics intrinsically involves the switching manifold. We construct explicit piecewise linear examples of both configurations. For crossing Bykov cycles, we prove the existence of hyperbolic invariant sets carrying symbolic dynamics, as in the classical smooth setting. For sliding Bykov cycles, under suitable open conditions, we obtain a countable symbolic dynamics on an invariant set; we prove that this invariant set has several interesting properties: the Hausdorff dimension of it and of its closure are equal, and it lies strictly between zero and one, and so both sets have zero Lebesgue measure; its closure is a Cantor set; and we also obtain conformal and invariant ergodic measures for the associated return dynamics. Finally, we show that convolution regularization of a suitable piecewise constant vector field produces a spatial cusp whose versal unfolding contains classical Bykov cycles. Thus, Bykov-type dynamics arises in the non-smooth setting both intrinsically, through crossing and sliding mechanisms, and through smooth regularization.
\end{abstract}

\section{Introduction}

Among the many heteroclinic configurations on $\R^3$ that have been studied, (smooth) Bykov cycles are of special interest; they consist of two hyperbolic saddle-focus equilibria with complementary Morse indices (the dimension of the unstable manifold) connected by two heteroclinic connections, with the intersection of the $2$-dimensional invariant manifolds being transverse. This structure generates an interaction between the local spiraling dynamics and the global heteroclinic structure, creating a rich dynamics.

However, the smooth theory of differential equations may not be enough to understand many mathematical models in contexts such as mechanical systems with impact, relay control systems, electrical switching circuits and biological models, where threshold effects are present. For these, the study of piecewise smooth vector fields is best suited.

In this paper we are going to consider mainly two questions regarding the relationship between Bykov cycles and piecewise smooth vector fields. For that, we will study mainly Filippov fields of the form
\begin{equation}\label{eq:filippov_system}
	Z(p) := 	\begin{cases}
			X(p), \qquad\text{if } f(p) > 0; \\
			Y(p), \qquad\text{if } f(p) < 0; 
			\end{cases}
\end{equation}
where $X$ and $Y$ are smooth vector fields, $f$ is a smooth function with $0$ as a regular value, and
$
\Sigma := f^{-1}(0)
$
denotes the switching set separating the two smooth regions. Other configurations will be also considered, and will be introduced as needed.

In the smooth theory, the trajectories evolve under a single differential equation. This is not the case for Filippov systems, which may exhibit several distinct regimes in the switching set; in particular, trajectories may cross, slide along or escape from it. For the case of sliding dynamics on $\Sigma$, this is completely determined by the Filippov sliding vector field, the convex combination of $X$ and $Y$ that is tangent to $\Sigma$. This generates an additional geometrical component for the dynamics that interacts with the two open smooth regimes in \eqref{eq:filippov_system}.

In particular, it can create new configurations of homoclinic and heteroclinic structures, some of them being objects inherent to this sliding dynamics, while also sharing much of the structure of some smooth counterparts; for instance, the Shilnikov cycle, a homoclinic saddle--focus structure, had a version inherent to Filippov system introduced by \cite{NovaesTeixeira2019ShilnikovProblemFilippov} and studied more deeply in \cite{NovaesPonceVarao2017ChaosInducedSliding, CunhaNovaesPonce2024HausdorffDimensionCantor}, called \emph{sliding Shilnikov cycles}.  Motivated by this, we are led to the natural question of whether and how the notion of Bykov cycles can be extended to the context of piecewise smooth dynamical systems. As far as we know, no such attempts to study Bykov cycles related to Filippov systems have been made so far.

This paper suggests two different configurations for this extension, according to the interaction between the heteroclinic trajectories and the switching set. The first type, called \emph{crossing Bykov cycles}, corresponds to the occurrence of both heteroclinic connections intersecting $\Sigma$ exclusively on crossing regions. Geometrically these are close to the classical smooth context, since the switching set acts only as an interface where the trajectories pass continuously. This case is illustrated in Figure~\ref{fig:crossing_bykov}. The second type, called \emph{sliding Bykov cycles}, is genuinely a piecewise smooth phenomenon, since it is given by the sliding dynamics on $\Sigma$. This occurs because in this case one of the trajectories of the heteroclinic cycle evolves along the Filippov sliding vector field before leaving the switching set through a visible fold point. This case is illustrated in Figure~\ref{fig:sliding_bykov}. We emphasize that the distinction between crossing and sliding Bykov cycles is not merely geometric, and that these two configurations generate essentially different dynamical mechanisms and require substantially different techniques and approaches for their analysis.

We also analyze the interplay between piecewise smooth vector fields, regularizations and the emergence of smooth Bykov cycles.

The paper is organized as follows: Section~\ref{section:preliminaries} presents the notion of the classical smooth Bykov cycles, then recalls the basic notions of Filippov systems in order to introduce the definitions of crossing and sliding Bykov cycles, and is where the main results are stated. Section~\ref{section:crossing_bykov_proof} investigates crossing Bykov cycles and establishes the
existence of Smale horseshoes and Section~\ref{section:crossing_bykov_examples} gives a family of examples of crossing Bykov cycles. Meanwhile, Section~\ref{section:sliding_bykov_proof} is devoted to sliding Bykov cycles and develops the symbolic
dynamics generated by the associated countable iterated function system in order to deduce geometric and dynamical properties of the system, and Section~\ref{section:sliding_bykov_examples} gives an example of a system containing a sliding Bykov cycle.
Then Section~\ref{section:foliations} studies triple switching configurations, convolution regularization and the resulting spatial cusp singularities together with
their versal unfoldings. Finally, in Section~\ref{section:conclusions} there are some concluding remarks of the investigation. We also discuss in Appendix~\ref{appendix:asymptotic_of_remainder} some results regarding estimates for remainder terms utilized throughout the proof; in Appendix~\ref{appendix:transversality} we address the transversality between the flow of the sliding vector field in $\Sigma$ and the projected curve, a critical component whose inclusion in the main text would disrupt the logical progression; in Appendix~\ref{appendix:cifs}, we give details of the relevant aspects of the theory of CIFS.

\begin{figure}[h!]
\centering
\begin{overpic}[width=0.5\linewidth]{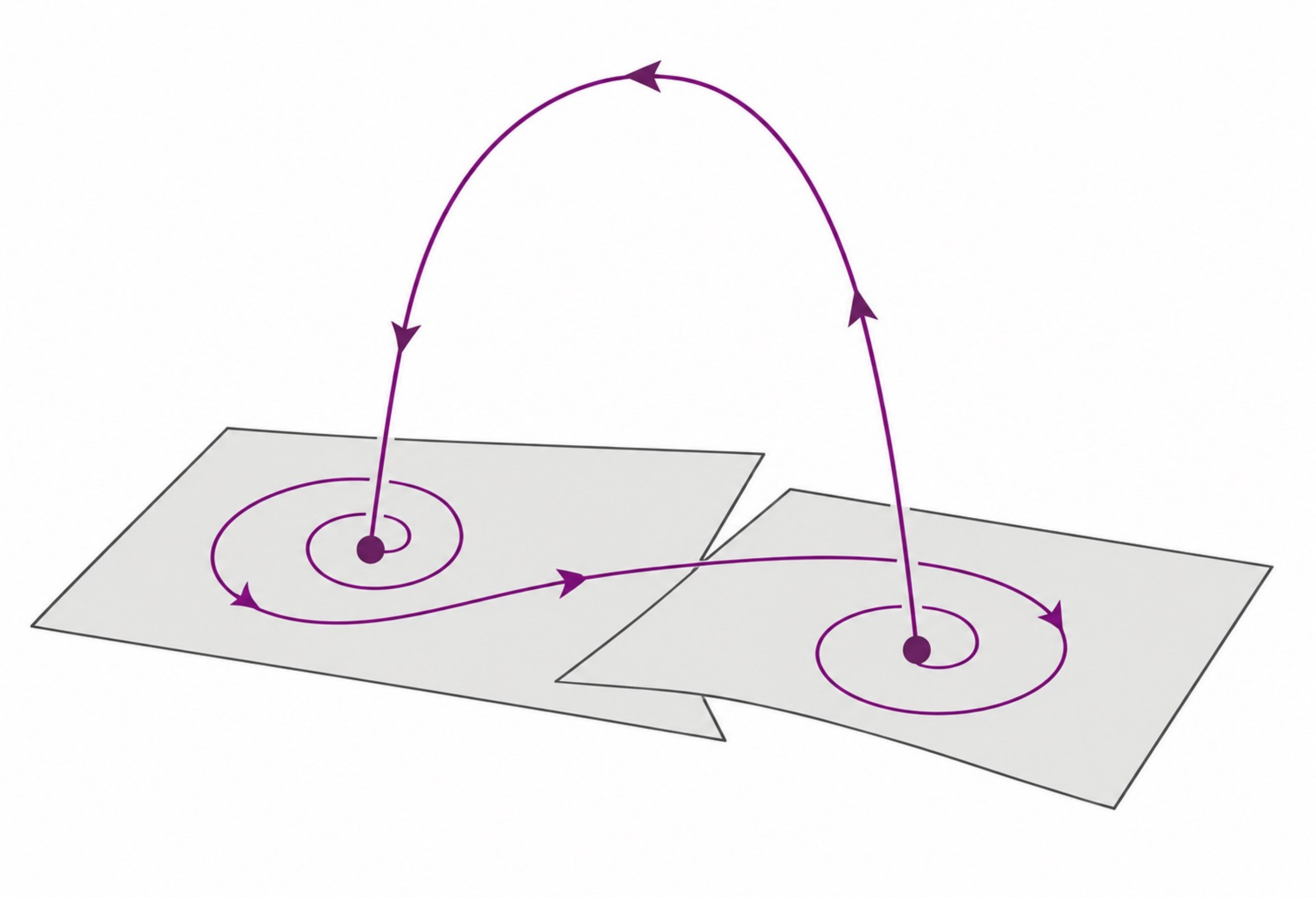}
\put(5,15){$W^u(p_u)$}
\put(90,8){$W^s(p_s)$}
\put(14,40){$W^s(p_u)$}
\put(70,35){$W^u(p_s)$}
\put(50,65){$\Gamma_{p_s\to p_u}$}
\put(40,27){$\Gamma_{p_u\to p_s}$}
\put(70,15){$p_s$}
\put(29.5,28){$p_u$}
\end{overpic}
\caption{A graphical representation of a (smooth) Bykov cycle.}
\label{fig:crossing_sliding_escaping}
\end{figure}

\section{Preliminaries and Main Results}\label{section:preliminaries}

In this section we introduce the necessary ingredients for defining Bykov cycles, as well as those for Filippov systems. With these, we can then discuss the notion of Bykov cycles in Filippov systems, allowing us to define the two structures of crossing and sliding Bykov cycles. Finally, we state our main results.

\subsection{Bykov Cycles} 
Our main objects of study are \emph{Bykov cycles}, for which we now give the definition for the classical smooth case.

\begin{definition}[Bykov Cycles]
Consider a smooth vector field $Z \in C^1(\R^3)$ with equilibrium points $p_s$ and $p_u$. Also, there exist open curves $\Gamma_{p_u\to p_s}$ connecting $p_u$ to $p_s$, and $\Gamma_{p_s\to p_u}$ connecting $p_s$ to $p_u$, trajectories of $Z$ with the following properties.

\begin{enumerate}[label = (B\arabic*), ref = (B\arabic*)]

\item\label{B1}
The point $p_s$ is a hyperbolic saddle-focus with eigenvalues of $DZ(p_s)$ given by
$\lambda_s$ and $-a_s \pm i b_s$ and the point $p_u$ is a hyperbolic saddle-focus with eigenvalues of $DZ(p_u)$ given by
$-\lambda_u$ and $a_u \pm i b_u$, where $\lambda_j, a_j > 0$ for $j\in\set{u,s}$, and $b_u b_s\neq 0$;

\item\label{B2}
The connection $\Gamma_{p_u\to p_s}$ satisfies
$
\Gamma_{p_u\to p_s} \subset W^u(p_u) \pitchfork W^s(p_s),
$
where $\pitchfork$ means that the intersection is transversal;

\item\label{B3}
The connection $\Gamma_{p_s\to p_u}$ satisfies
$
\Gamma_{p_s\to p_u} \subset W^u(p_s)\cap W^s(p_u).
$
\end{enumerate}

Then the cycle $\Gamma := \{p_u\} \cup \Gamma_{p_u\to p_s} \cup \{p_s\} \cup \Gamma_{p_s\to p_u}$ is called a \emph{Bykov cycle}.
\end{definition}

We remark that  the signs of $b_s$ and $b_u$ are fixed according to the chosen orientations of the corresponding two-dimensional eigenspaces, with $b_s$ and $b_u$ denoting the corresponding angular velocities.

For the next result, we denote by $\mathcal{B}\subset C^2(\R^3)$ the class of vector fields possessing
Bykov cycles. One of the main properties of Bykov cycles is the following.
\begin{theorem}\label{ppst:smooth_bykov_cycles}
There exists an open set $\mathcal{E}\subset\mathcal{B}$ such that, for every $Z\in\mathcal{E}$ and every tubular neighborhood $U$ of its Bykov cycle $\Gamma$ contains a cross-section $S = S(U)\subset U$ and a sequence of hyperbolic compact sets
$\{\Lambda_k\}_{k\ge2}$, with $\Lambda_k\subset S$, such that
$
\pi|_{\Lambda_k}\colon\Lambda_k\rightarrow\Lambda_k
$
is topologically conjugate to the shift on $k$ symbols.
\end{theorem}

The proof can be found in \cite[Theorem~3.2]{BarrientosIbanezRodrigues2019EmergenceChaoticDynamics} and \cite[Theorem~2 and Section~5.3]{LabouriauRodrigues2015DenseHeteroclinicTangencies}.

The existence and behavior of heteroclinic cycles are of great interest in dynamical systems, making a rigorous understanding of Bykov cycles essential. With this in mind, we now investigate the existence of analogous cycles in piecewise smooth settings.

\subsection{Piecewise Smooth Vector Fields in $\R^3$}

A \emph{piecewise smooth vector field} in $\R^3$ is a vector field defined by 
smooth expressions on a finite partition of the phase space, separated by a switching set.

Let $\mathcal{F} = \set{f_i:\R^3 \to \R}_{1\le i\le k}$ be a collection of smooth functions and define the \emph{switching set} as
$
\Sigma := \bigcup_{i=1}^k \Sigma_i,
$
with $\Sigma_i := f_i^{-1}(0)$ and where $0$ is a regular value for all functions in $\cal{F}$, that is,
$
\nabla f_i(p) \neq 0 \qquad \text{for all } p \in \Sigma_i.
$

The complement $\R^3 \setminus \Sigma$ is partitioned into finitely many regions
$
\{R_\alpha\}_{\alpha \in A},
$
where each region is determined by the signs of the functions $f_i$. We will call the smooth regions
$
R := \R^3\setminus\Sigma = \cup_{\alpha\in A}R_{\alpha}.
$

A piecewise smooth vector field is a map $Z:\R^3 \to \R^3$ such that
$
Z|_{R_\alpha} = Z_\alpha,
$
where each $Z_\alpha$ is of class $C^{l}$ on a neighborhood of $\cl{R_\alpha}$, with $l\geq 1$ an integer.
The set of piecewise smooth vector fields induced by $\cal{F}$ in this setting is denoted by $\Omega_{\mathcal{F}}^{l}$.
It has a natural inclusion in $(C^{l})^{\abs{A}}$, allowing us to endow it with the product topology; the case where $\cal{F} = \emptyset$ is just the classical smooth case.

Let us assume the particular case where $\Sigma = f^{-1}(0)$, where $0$ is a regular value for a smooth function $f:\R^3 \to \R$. This induces the classical two-zone setting
$
R^+ = \{f>0\}$ and $R^- = \{f<0\},
$
where the vector field can be written as
\[
Z =
\frac{1+\sgn(f)}{2}\, X
+
\frac{1-\sgn(f)}{2}\, Y.
\]

While we may have $\Sigma$ being formed by many functions $\set{f_i}$, whenever  $p\in\Sigma_i\subset\Sigma$, with $p\not\in\Sigma_i\cap\Sigma_j$ for $i\neq j$ (which are called \emph{generic} points), we have locally the two-zone setting and we will denote these regions and vector fields by $\Sigma_i$, $R^+_i$, $R^-_i$, $X_i$, $Y_i$ and $Z_i$, but if there is no risk of confusion we will not use the index.

\subsubsection{Regular Switching Points and Switching Regions}

Assume that the switching point is generic. A point $p \in \Sigma$ is called a \emph{regular switching point} if
$
F f(p) \neq 0,
$ where $
F f(p) = \inner{F(p), \nabla f(p)}.
$

The dynamics near such points is classified as
\begin{itemize}
\item $\Sigma^c := \{p\in\Sigma \text{ such that } Xf(p)\cdot Yf(p) > 0\}$;
\item $\Sigma^s := \{p\in\Sigma \text{ such that } Xf(p) < 0$ and $Yf(p) > 0\}$;
\item $\Sigma^e := \{p\in\Sigma \text{ such that } Xf(p) > 0$ and $Yf(p) < 0\}$;
\item $\Sigma^t := \{p\in\Sigma \text{ such that } Xf(p)\cdot Yf(p) = 0\}$.
\end{itemize}
The region $\Sigma^c$ is called the \emph{crossing region}, the region $\Sigma^s$ is called the \emph{sliding region}, the region $\Sigma^e$ is called the \emph{escaping region} and the region $\Sigma^t$ is called the \emph{tangent region}.

\begin{figure}[h!]
\centering
\begin{overpic}[width=0.8\linewidth]{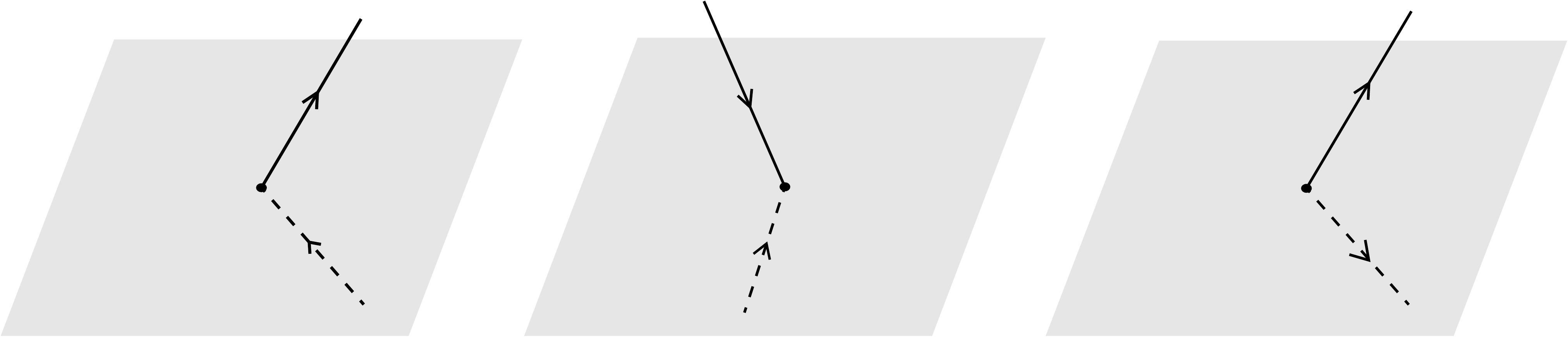}
\put(14,9){$p$}
\put(47,9){$p$}
\put(80.5,9){$p$}
\put(4,1){$\Sigma^c$}
\put(38,1){$\Sigma^s$}
\put(71,1){$\Sigma^e$}
\end{overpic}
\caption{A graphical representation of the crossing, sliding and escaping regions on the switching manifold.}
\label{fig:crossing_sliding_escaping}
\end{figure}
%

\subsubsection{Fold Points}

A point $p \in \Sigma$ is a \emph{regular fold point} of $X$ if
$
Xf(p) = 0
$
and
$
X^2f(p) \neq 0,
$
where
$
X^2f(p) = X(Xf)(p).
$
It is a \emph{visible fold-regular point} if $X^2f(p) > 0$ and $Yf(p) \neq 0$ (note that in the boundary of sliding regions, we must have $Yf(p) > 0$ for this to occur). All these definitions can be made analogous for $Y$.

\subsubsection{Sliding Vector Field}

On the sliding and escaping regions the points are generic, so in a small neighborhood they are uniquely determined as contained in $\Sigma_i$, for some $i$; the induced local dynamics on the switching set $\Sigma$ are governed by the sliding vector field
\[
Z_{\Sigma}(p) = Z_{{\Sigma},i}(p) :=
\frac{Y_if_i(p)\, X_i(p) - X_if_i(p)\, Y_i(p)}
{Y_if_i(p) - X_if_i(p)}.
\]

Because this vector field is strictly tangent to $\Sigma$, it naturally defines the dynamics within the sliding and escaping regions. This leads to a specific type of equilibrium for the Filippov system $Z$, given by the equilibria of $Z_{\Sigma}$. These are points $p^*\in \Sigma^{s,e}$ where $Z_{\Sigma}(p^*)=0$, known as \emph{pseudo--equilibria of $Z$}.

A pseudo--equilibrium is said to be \emph{hyperbolic} if it is a hyperbolic equilibrium of the planar vector field $Z_{\Sigma}$. Furthermore, if $p^* \in \Sigma^s$ is an unstable hyperbolic focus of $Z_{\Sigma}$, or if $p^* \in \Sigma^e$ is a stable hyperbolic focus of $Z_{\Sigma}$, then $p^*$ is defined as a \emph{hyperbolic pseudo-saddle-focus}.

\subsection{Bykov Cycles in Piecewise Smooth Systems}
Having introduced both the notion of smooth Bykov cycles and piecewise smooth vector fields, we are now in a position to provide the definition of \emph{piecewise smooth Bykov cycles} in different configurations, which expands the notion for the piecewise smooth context.

\begin{figure}[h!]
\centering
\begin{subfigure}{.5\textwidth}
  \centering
  	\begin{overpic}[width=0.6\linewidth]{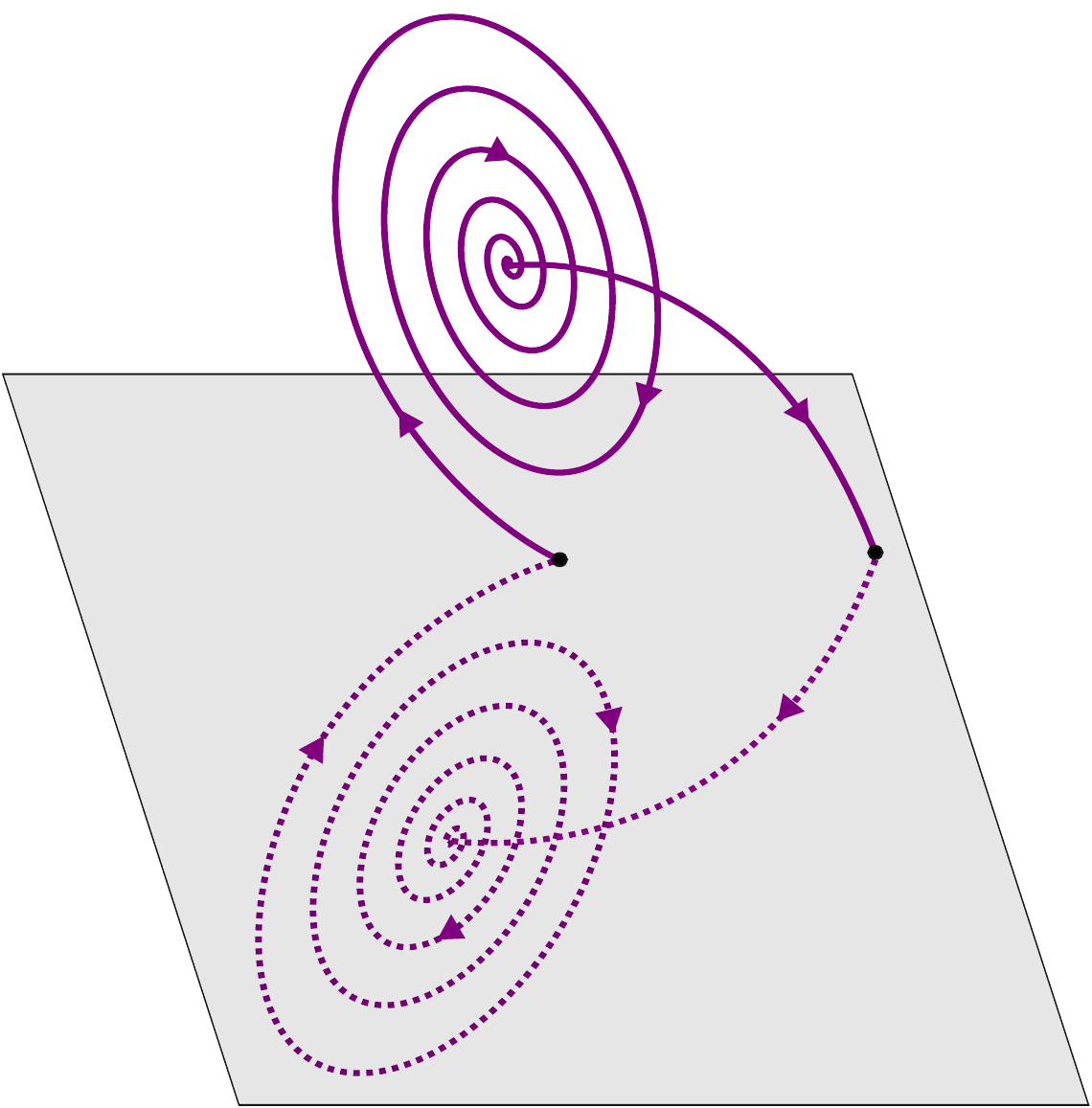}
	\put(7,58){$\Sigma^c$}
	\put(56,90){$\Gamma_{p\to p_s}$}
	\put(68,70){$\Gamma_{p_s\to q}$}
	\put(52,10){$\Gamma_{p_u\to p}$}
	\put(68,25){$\Gamma_{q\to p_u}$}
	\put(52,48){$p$}
	\put(73,48){$q$}
	\end{overpic}
  \caption{A crossing Bykov cycle.}
  \label{fig:crossing_bykov}
\end{subfigure}%
\begin{subfigure}{.5\textwidth}
  \centering
	\begin{overpic}[width=0.8\linewidth]{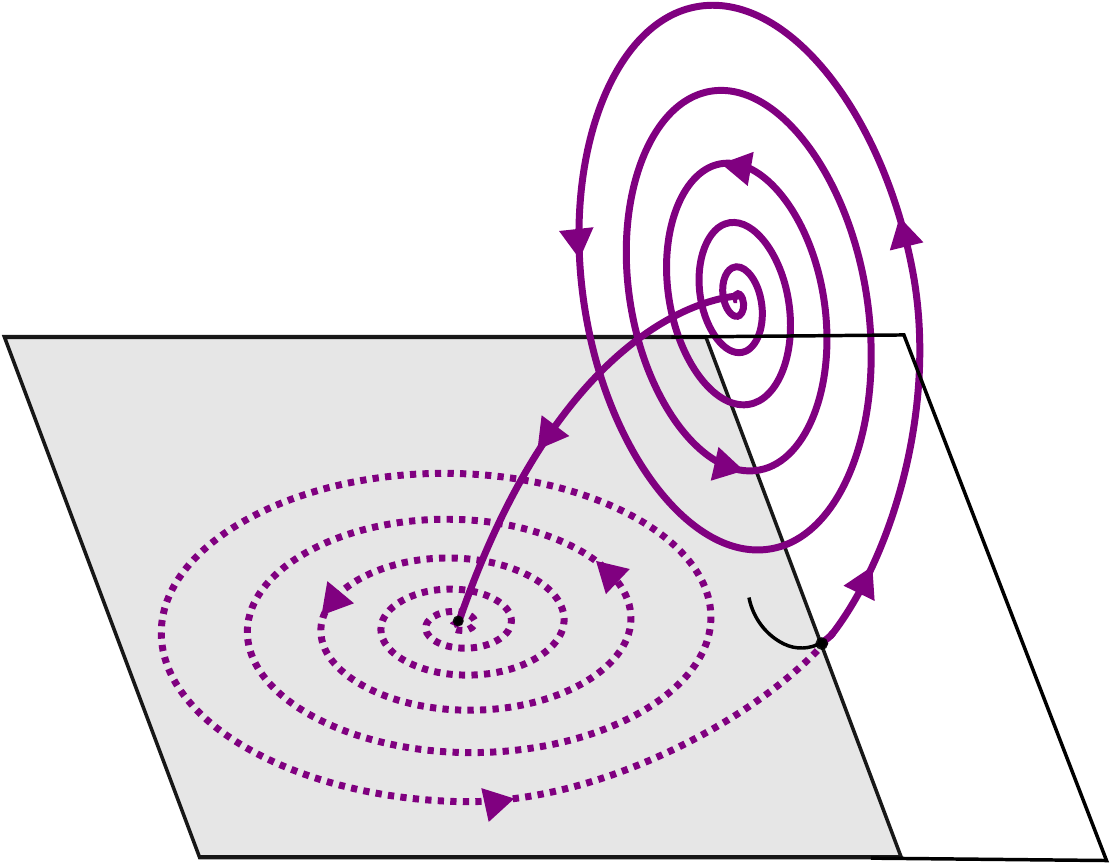}
	\put(7,40){$\Sigma^s$}
	\put(76,17){$q$}
	\put(23,4){$\Gamma_{p_u\to q}$}
	\put(85,50){$\Gamma_{q\to p_s}$}
	\put(33,40){$\Gamma_{p_s\to p_u}$}
	\put(69,3){$\partial\Sigma^s$}
	\end{overpic}
  \caption{A sliding Bykov cycle.}
  \label{fig:sliding_bykov}
\end{subfigure}%
\caption{The two configurations of piecewise smooth Bykov cycles.}
\label{fig:piecewise_smooth_bykov}
\end{figure}

For that, given a sufficiently small neighborhood $U\subset\Sigma$ of $p\in\Sigma$, we define the sets
\begin{equation*}
W_\Sigma^s(p) := \set{x\in U : \flow_X^t(x) \to p \text{ as } t\to +\infty}
\end{equation*}
and
\begin{equation*}
W_\Sigma^u(p) := \set{x\in U : \flow_Y^t(x) \to p \text{ as } t\to -\infty}.
\end{equation*}

Now we can make the following definition.

\begin{definition}[Crossing Bykov Cycles]\label{def:crossing_bykov}
Let $Z = (X,Y) \in \Omega^1(\R^3)$ be a Filippov system with switching set $\Sigma = f^{-1}(0)$. Suppose that it has two equilibrium points $p_s\in R^+$ and $p_u\in R^-$. Also, there exist open curves $\Gamma_{p_u\to p_s}$ connecting $p_u$ to $p_s$, and $\Gamma_{p_s\to p_u}$ connecting $p_s$ to $p_u$, trajectories of $Z$ with the following properties.


\begin{enumerate}[label = (CB\arabic*), ref = (CB\arabic*)]

\item\label{CB1}
The point $p_s$ is a hyperbolic saddle-focus for $X$ with eigenvalues of $DX(p_s)$ given by
$\lambda_s$ and $-a_s \pm i b_s$, and the point $p_u$ is a hyperbolic saddle-focus for $Y$ with eigenvalues of $DY(p_u)$ given by $-\lambda_u$ and $a_u \pm i b_u$, where $\lambda_j, a_j > 0$ for $j\in\set{u,s}$, and $b_u b_s\neq 0$;


\item\label{CB2}
The connection $\Gamma_{p_u\to p_s}$ can be decomposed as $\Gamma_{p_u\to p_s} = \Gamma_{p_u\to p} \,\dot\cup\, \set{p} \,\dot\cup\, \Gamma_{p\to p_s}$,
with
$\Gamma_{p_u\to p} \subset R^-$, $p\in\Sigma^c$ and $\Gamma_{p\to p_s}\subset R^+$.

\noindent Besides, $W_\Sigma^s(p_s)\pitchfork W_\Sigma^u(p_u)\subset\Sigma^c$, where $\pitchfork$ means that the intersection is transversal;


\item\label{CB3}
The connection $\Gamma_{p_s\to p_u}$ can be decomposed as
$
\Gamma_{p_s\to p_u} = \Gamma_{p_s\to q} \,\dot\cup\, \set{q} \,\dot\cup\, \Gamma_{q\to p_u}
$
with
$\Gamma_{p_s\to q} \subset R^+$, $q\in\Sigma^c$ and $\Gamma_{q\to p_u}\subset R^-$.
\end{enumerate}


Then the cycle $\Gamma := \{p_u\} \cup \Gamma_{p_u\to p_s} \cup \{p_s\} \cup \Gamma_{p_s\to p_u}$ is called a \emph{crossing Bykov cycle}.
\end{definition}

A picture illustrating a crossing Bykov cycle can be found in Figure~\ref{fig:crossing_bykov}.

There is another object that shares some of the characteristics of Bykov cycles, while having dynamics intrinsically linked to Filippov systems. We now give the definition.

\begin{definition}[Sliding Bykov Cycles]
\label{def:sliding_bykov}
Let $Z = (X,Y) \in \Omega^1(\R^3)$ be a Filippov system with switching set $\Sigma = f^{-1}(0)$. Suppose that it has an equilibrium $p_s\in R^+$ and a pseudo-equilibrium $p_u\in \Sigma^s$. Also, there exist open curves $\Gamma_{p_u\to p_s}$ connecting $p_u$ to $p_s$, and $\Gamma_{p_s\to p_u}$ connecting $p_s$ to $p_u$, which are trajectories of $Z$, with the following properties.
\begin{enumerate}[label = (SB\arabic*), ref = (SB\arabic*)]

\item\label{SB1}
The point $p_s$ is a hyperbolic saddle-focus for $X$ with eigenvalues of $DX(p_s)$
$
\lambda_s$ and $-a_s \pm i b_s,
$
and the point $p_u$ is a hyperbolic pseudo-saddle-focus for $Z_{\Sigma}$ with eigenvalues of $DZ_{\Sigma}(p_u)$
$
a_u \pm i b_u,
$
where $\lambda_s, a_j > 0$ for $j\in\set{u,s}$, $b_u b_s\neq 0$ (here the signs of $b_s$ and $b_u$ are fixed according to the chosen orientations of the corresponding two-dimensional eigenspaces, with
$b_s$ and $b_u$ denoting the corresponding angular velocities);


\item\label{SB2}
The connection $\Gamma_{p_u\to p_s}$ can be decomposed as
$
\Gamma_{p_u\to p_s} = \Gamma_{p_u\to q} \,\dot\cup\, \set{q} \,\dot\cup\, \Gamma_{q\to p_s}
$
(each one a connected component), with $\Gamma_{p_u\to q} \subset \Sigma^s$, $q\in\partial\Sigma^s$ and $\Gamma_{q\to p_s}\subset R^+$.

\noindent Besides, $q$ is a visible fold-regular point, being a visible fold for $X$;


\item\label{SB3}
The connection $\Gamma_{p_s\to p_u}$ is a trajectory of $X$ satisfying
$
\Gamma_{p_s\to p_u} \subset W^u(p_s)\cap R^+,
$
such that its forward closure intersects the switching set exactly at the pseudo-equilibrium, that is, $\overline{\Gamma_{p_s\to p_u}} \cap \Sigma = \set{p_u}$.
\end{enumerate}


Then the cycle $\Gamma := \set{p_u} \cup \Gamma_{p_u\to p_s} \cup \set{p_s} \cup \Gamma_{p_s\to p_u}$ is called a \emph{sliding Bykov cycle}.
\end{definition}

A picture illustrating a sliding Bykov cycle can be found in Figure~\ref{fig:sliding_bykov}.

In what follows, we state the main results about these objects.

\subsection{Main Results}

As usual in the understanding of dynamical systems, the precise statements refer to properties of a first-return map; so let $S$ be a set and consider a map $\pi:\Dom{\pi}\subset S\to S$, then for any non-empty subset $U\subset S$, we consider the \emph{local invariant map} and the \emph{local invariant set} given by
\begin{equation}\label{eq:local}
	\pi_U := \pi|_{\Dom{\pi}\cap U}
\qquad\text{ and }\qquad
\Lambda_U := \set{w\in U : \pi^k(w)\in \Dom{\pi}\cap U, \text{for all }k\ge0}.
\end{equation}

The first result concerns the crossing Bykov cycles, stating that its resemblance to the smooth case induces similar dynamical properties between them. For that, we consider $\mathcal{B}\subset\Omega^2(\R^3)$ as the class of Filippov systems possessing crossing Bykov cycles as given in Definition~\ref{def:crossing_bykov}, and for $Z\in\mathcal B$ we denote $\Gamma$ as its corresponding cycle. We investigate the dynamics of the first-return map near $\Gamma$.

\begin{maintheorem}\label{thm:crossing_bykov}
For every $Z\in\mathcal{B}$ and every tubular neighborhood $U$ of its crossing Bykov cycle $\Gamma$, there exists
a cross-section $S = S(U)\subset U$ of $\Gamma$ and a sequence of hyperbolic compact sets
$\set{\Lambda_k}_{k\ge2}$, with $\Lambda_k\subset S$, such that
$
\pi|_{\Lambda_k}\colon\Lambda_k\rightarrow\Lambda_k
$
is topologically conjugate to the shift on $k$ symbols.
\end{maintheorem}

The next result focuses on sliding Bykov cycles. We restrict our focus to Filippov systems $Z = (X,Y)$ with a condition on the local linearizations.

\begin{enumerate}[label=(P\arabic*), ref=(P\arabic*), start=1]
\item\label{cond:linearization} The vector field $X$ at $p_s$ and the sliding vector field $Z_\Sigma$ at $p_u$ admit local $C^2$ conjugacies $h_s$ and $h_u$, respectively, with $C^2$ inverses, to their linear parts.
\end{enumerate}

This may be achieved, for instance, if $X$ and $Z_\Sigma$ are
$C^\infty$ in neighborhoods of $p_s$ and $p_u$, respectively, and
$\lambda_s\neq2a_s$, by an application of Samovol's Linearization Theorem (see \cite{Samovol1988SmoothLinearizationSystems}).

We consider a sufficiently small $r>0$ and the set
\begin{equation}\label{eq:gamma}
\gamma = \gamma_r := \cl{B_r(q)\cap\partial\Sigma^s},
\end{equation}
which is a compact curve containing the visible fold-regular point $q\in\partial\Sigma^s$.

Furthermore, we discuss some quantities referring to a certain map $\hat{\mu}_{glob}$ and a certain point $\hat{u}_{out}$, both related to the connection $\Gamma_{p_s\to p_u}$. This map is defined and studied in detail in Section~\ref{section:sliding_bykov_proof}. We then denote
\begin{equation}\label{eq:L_quantities}
L := D\hat{\mu}_{glob}(\hat{u}_{out}),
\qquad
S := \sgn(\det L)
\qquad\text{ and }\qquad
M := \frac{\sigma_1^2 + \sigma_2^2}{2 \sigma_1 \sigma_2} \ge 1,
\end{equation}
where $\sigma_1 \ge \sigma_2 > 0$ are the singular values of $L$ (which can be characterized as
$
\sigma_1=\max_{\|v\|=1}\|Lv\|
$
and
$
\sigma_2=\min_{\|v\|=1}\|Lv\|
$), and $M\ge1$ is a measure of the distortion of the ellipse induced by the linear map $L$. Furthermore, we define
$
c_s := (-a_s)/b_s$ and $c_u := a_u/b_u,
$
and then state the conditions
\begin{enumerate}[label=(P\arabic*), ref=(P\arabic*), start=2]
	\item\label{cond:non_resonance_undistorted}
	$c_s - Sc_u \neq 0$;
    \item\label{cond:distortion_bounds}
    $\abs{M - S c_s c_u} < \sqrt{(c_s^2 + 1)(c_u^2 + 1)}$.
\end{enumerate}

Here, \ref{cond:non_resonance_undistorted} is a non-resonance condition, while \ref{cond:distortion_bounds} is used to ensure that the projected spiral in $\Sigma$ is transversal to the flow of the sliding vector field near the pseudo-saddle-focus.

In order to state the main result for sliding Bykov cycles, we now consider a Filippov system $Z=(X,Y)\in\Omega^3(\mathbb R^3)$ possessing a sliding Bykov cycle $\Gamma$ with distinguished point $q$, as described in Definition~\ref{def:sliding_bykov}, and satisfying conditions \ref{cond:linearization}--\ref{cond:distortion_bounds}. Let $\pi$ be the associated first-return map on $\gamma$, as defined in \eqref{eq:gamma}, and for every sufficiently small neighborhood $U\subseteq\gamma$ of $q$, we consider the corresponding local return map $\pi_U$ and its invariant set $\Lambda_U$, defined as in \eqref{eq:local}.
We denote $\hausdorff{\cdot}$ for the Hausdorff dimension and $\opname{m}_1(\cdot)$ for the one-dimensional Lebesgue measure. The definition of conformal measures in a more detailed context can be found in Definition~\ref{definition:t_conformal_measures}. In these conditions, we have the following results.

\begin{maintheorem}\label{thm:sliding_bykov}
For every sufficiently small neighborhood $U\subset\gamma$ of $q$, the local map $\pi_U$ and the local invariant set $\Lambda_U$ satisfy:
\begin{enumerate}[label= (\alph*), ref = \ref{thm:sliding_bykov}(\alph*)]
	
	\item\label{thm:symbolic_dynamics} \emph{Symbolic Dynamics}: $\pi_U|_{\Lambda_U}:\Lambda_U\to\Lambda_U$ is topologically conjugate to the full shift on a countable alphabet $(\sigma, \N^\N)$.
	\item\label{thm:fractal_geometry} \emph{Set and Fractal Structure}: $0<\hausdorff{\Lambda_U} = \hausdorff{\ov\Lambda_U} < 1$ and $\opname{m}_1(\Lambda_U) = \opname{m}_1(\ov\Lambda_U) =0$; besides,  $\overline{\Lambda_U} = \Lambda_U\,\dot\cup\,Q_U$, where $Q_U := \cup_{k\ge0}\pi_U^{-k}(q)$ is a countable set, and it is a Cantor set.
	\item\label{thm:ergodic_properties} \emph{Ergodic Properties}: There exists a probability measure $m_U$ on $\Lambda_U$ that is $t_U$-conformal
	\footnote{A probability measure $m$ on a set $\Delta\subset K$ is called a \emph{$t$-conformal measure} for a countable family of contractions $\set{f_i:K\to K}_{i\in I}$ if it satisfies the following properties:
	\begin{enumerate}
    \item $m(\Delta) = 1$;
    \item $m(f_i(B)) = \int_{B}\norm{Df_i}^tdm$, for all $i\in I$ and for all Borel sets $B\subset\Delta$;
    \item $m(f_i(K)\cap f_j(K)) = 0$, for all $i\in I$, $i\neq j$.
	\end{enumerate}
	}
	for the countable family of inverse branches of $\pi_U$, for $t_U := \hausdorff{\Lambda_U}$. Moreover, there exists a unique measure $m_U^* \sim m_U$ which is invariant and ergodic under $\pi_U|_{\Lambda_U}$.
	
	\end{enumerate}
\end{maintheorem}

We remark that the techniques used in Theorem~\ref{thm:fractal_geometry} allow us to produce sharper estimates for the Hausdorff dimension. This will be discussed in Remark~\ref{rem:sharper_bounds}.

There is yet another way that Filippov systems and Bykov cycles are related, which is the way the latter may appear when trying to regularize the former. Here we introduce the notion, leaving the details for Section~\ref{section:foliations}.

\emph{Regularizations} of Filippov systems are techniques that associate a family of smooth systems that approximate the original piecewise smooth vector field in some sense. Subsection~\ref{subsection:convolution_regularization} presents one of these techniques, for instance.

Our third result concerns the regularization of piecewise constant
vector fields near a triple crossing of the switching set. More precisely,
we consider collections
\[
C=\{C_{\widetilde v}\}_{\widetilde v\in\{-1,1\}^3}
\]
of constant vector fields, one for each region determined by the three
switching hypersurfaces, and denote by $Y_C$ the smooth vector field
obtained from the corresponding discontinuous system through convolution
regularization.

Although the vector fields defining the original piecewise system are
constant, their regularization may exhibit nontrivial local dynamics near
the triple crossing. The following result shows that the collection $C$
can be chosen so that $Y_C$ has a spatial cusp at the origin. Moreover,
perturbations of $C$ generate a versal unfolding of this singularity,
within which a classical Bykov cycle occurs.

\begin{maintheorem}\label{thm:versal_unfolding}
There exists a collection $C=\{C_{\tilde{v}}\}$ of constant vector fields such
 that the corresponding vector field $Y_C$ has a spatial cusp at the origin. Moreover, 
 this spatial cusp unfolds versally under small perturbations of the collection $C$. Besides, this versal unfolding has a Bykov cycle.
\end{maintheorem}

\section{Crossing Bykov Cycles: Proof of Theorem \ref{thm:crossing_bykov}}\label{section:crossing_bykov_proof}

We sketch the proof of Theorem~\ref{thm:crossing_bykov} in two steps, and invoke the main mechanisms for smooth systems from established references to complete the proof. The steps are the following.

\begin{enumerate}
    \item[\textit{Step 1.}]
    \textit{Smoothness of Global Maps.}
    We show that while the trajectories from a neighborhood of $p_u$ to a neighborhood of $p_s$ may not be smooth when crossing the switching set, the global maps between them are.

    \item[\textit{Step 2.}]
    \textit{Adaptation of the Smooth Case.}
    Since the global maps are smooth, we compose the local and global maps, resulting in a smooth first-return map. Then using the same construction of the smooth Bykov cycles mapping horizontal strips into vertical strips, we verify that they create a horseshoe phenomenon. 
\end{enumerate}

So now we proceed with the proof.

\begin{figure}[h!]
\centering
\begin{overpic}[width=0.8\linewidth]{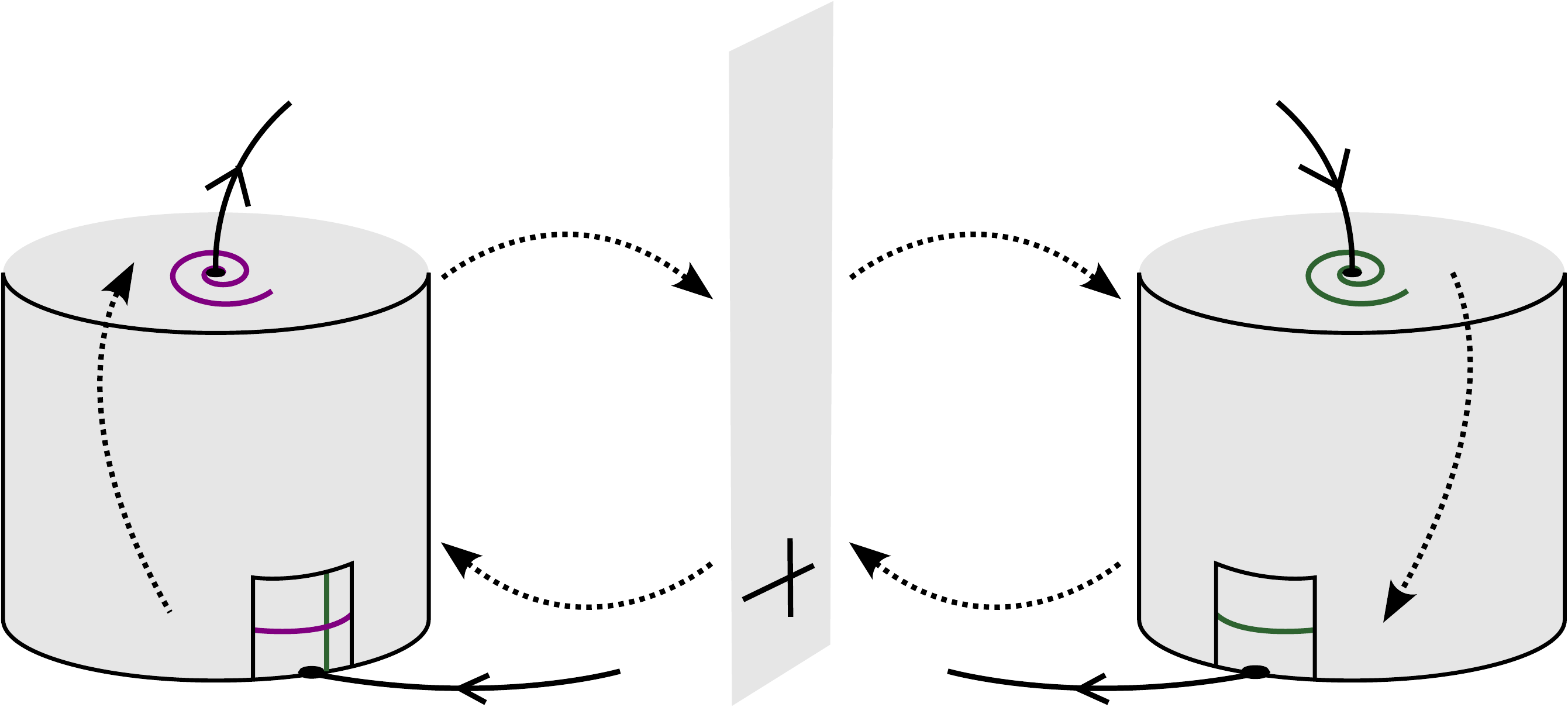}
\put(18,40){$\Gamma_{p_s\to p_u}$}
\put(74,40){$\Gamma_{p_s\to p_u}$}
\put(30,3){$\Gamma_{p_u\to p_s}$}
\put(65,3){$\Gamma_{p_u\to p_s}$}
\put(89,16){$\pi_u$}
\put(8,16){$\pi_s$}
\put(34,8){$\pi_{\Sigma\to s}$}
\put(61,8){$\pi_{u\to \Sigma}$}
\put(34,31.5){$\pi_{s\to \Sigma}$}
\put(61,31.5){$\pi_{\Sigma\to u}$}
\put(48,38){$\Sigma^c$}
\put(39.5,13){$W^u_{\Sigma}(p_u)\pitchfork W^s_{\Sigma}(p_s)$}
\end{overpic}
\caption{Graphical representation of the local and global maps of a crossing Bykov cycle. A horizontal strip in the wall of a cylindrical neighborhood of $p_s$ is mapped to a spiral in the top disk. It travels to the top disk of the cylindrical neighborhood of $p_u$, where it is mapped to a horizontal line again. Then, when traveling to the initial cylinder wall, the transversal intersection $W^u_{\Sigma}(p_u)\pitchfork W^s_{\Sigma}(p_s)$ generates the transversal encounter of the first-return map, which is the usual mechanism for the horseshoe phenomenon.}
\label{fig:horseshoe}
\end{figure}

\subsection{Smoothness of Global Maps}
Let $Z=(X,Y)\in\mathcal{B}\subset\Omega^2(\R^3)$ with a crossing Bykov cycle $\Gamma$. We consider a sufficiently small tubular neighborhood of $\Gamma$ and choose local cross-sections near $p_s$ and $p_u$.
The first-return map can be written as
\[
\pi
=
\pi_{u\to s}\circ\pi_u\circ
\pi_{s\to u}\circ\pi_s,
\]
where $\pi_s$ and $\pi_u$ are the local transition maps near
$p_s$ and $p_u$, respectively, while $\pi_{s\to u}$ and
$\pi_{u\to s}$ are the global transition maps along the two
heteroclinic connections.

First, we verify that the global transition maps have the
regularity required in the classical construction. For that, let us first consider the transition map $\pi_{s\to u}$ from $p_s$ to $p_u$. Since the connection crosses
$\Sigma$ at $q\in\Sigma^c$, we can write
\[
\pi_{s\to u}
=
\pi_{\Sigma\to u}\circ\pi_{s\to\Sigma},
\]
where $\operatorname{Im}(\pi_{s\to \Sigma}) = \operatorname{Dom}(\pi_{\Sigma\to u})\subset \Sigma$.
Let $w\in\Gamma_{p_s\to p_u}\cap R^-$ be a point in the connection. The time $t_Y(w)$ for the trajectory of $Y$
to reach $\Sigma$ is given implicitly by
$
f(\varphi_Y^{t_Y(w)}(w))=0.
$
At the image point $q\in\Sigma$, since it is a crossing region, we have
\[
\frac{\partial}{\partial t}
f(\varphi_Y^{t_Y(w)}(w))
=
\inner{
Y(q),
\nabla f(q)
}
\neq0.
\]
Hence, by the Implicit Function Theorem, $t_Y$ is of class $C^2$,
and therefore
$
\pi_{s\to\Sigma}(w)
=
\varphi_Y^{t_Y(w)}(w)
$
is a $C^2$ local diffeomorphism. Using the same argument, we can deduce that $\pi_{\Sigma\to u}$ is a $C^2$ local diffeomorphism, and consequently the composition $\pi_{s\to u}$ is a $C^2$ local diffeomorphism. Again, in the same way we can show that $\pi_{u\to s}$ is also a $C^2$ local diffeomorphism.

\subsection{Adaptation of the Smooth Case}
We have the following properties.
\begin{enumerate}
    \item the two local maps have the usual saddle-focus spiral geometry;
    \item the two global transition maps are local $C^2$ diffeomorphisms;
    \item by \textnormal{(CB2)}, the two-dimensional invariant
    manifolds meet transversely along
    $\Gamma_{p_u\to p_s}$;
\end{enumerate}
this allows us to apply to $\pi$ the same mechanism for the construction of Smale horseshoes of smooth Bykov cycles, which is essentially this: for a sufficiently small rectangle $R$ containing $\Gamma_{p_u\to p_s}$ in its border, the set of points returning to $R$ is contained in a sequence of horizontal strips $\set{h_{k}}_{k\in\N}$ that accumulates in $W^s_{loc}(p_s)$; each of these horizontal strips is locally mapped to a spiral near $\Gamma_{p_s\to p_u}$, which then returns and is mapped into horizontal strips; finally, the horizontal strips follow along $\Gamma_{p_u\to p_s}$, but since the unstable manifold of $p_u$ and the stable manifold of $p_s$ intersect transversally, they are mapped into vertical strips $\set{v_{k}}_{k\in\N}$ crossing transversally the original horizontal strips $\set{h_{k}}_{k\in\N}$, inducing a sequence of uniformly hyperbolic horseshoes accumulating on $\Gamma$.

In particular, for every $k\ge2$, we can choose $k$ sufficiently
close spiralling strips so that the associated invariant set
$\Lambda_k$ is a compact hyperbolic set and
$
\pi|_{\Lambda_k}:\Lambda_k\rightarrow\Lambda_k
$
is topologically conjugate to the shift on $k$ symbols.
Since the strips can be chosen arbitrarily close to the heteroclinic
cycle, the sets $\Lambda_k$ accumulate on $\Gamma\cap S\subset U$, for any
tubular neighborhood $U$ of $\Gamma$.
\qed

\bigskip
For the detailed proof, this method can be found in \cite[Theorem 3.2 and Section 3.4]{BarrientosIbanezRodrigues2019EmergenceChaoticDynamics}, \cite[Theorem 2 and Section 5.3]{LabouriauRodrigues2015DenseHeteroclinicTangencies}, and \cite[Corollary 10 and Section 5]{Rodrigues2013RepellingDynamicsBykov}.

\section{Crossing Bykov Cycles: Examples}\label{section:crossing_bykov_examples}
Now we provide a family of examples of piecewise smooth vector fields, even more piecewise linear vector fields, possessing a crossing Bykov cycle.

\begin{example}[Piecewise Linear Crossing Bykov Cycle]
\label{example:crossing_bykov}
We construct an explicit piecewise smooth vector field $Z = (X,Y) \in \Omega^2(\R^3)$ with switching set $\Sigma = f^{-1}(0)$, where $f(x,y,z) = z$, defined by:
\begin{equation}\label{eq:crossing_bykov}
    \begin{cases} 
        X(p) = A_s (p - p_s) & \text{for } p \in R^+ \ (z > 0), \\
        Y(p) = A_u (p - p_u) & \text{for } p \in R^- \ (z < 0),
    \end{cases}
\end{equation}
where $p_s = (0,0,1)^T$ and $p_u = (0,0,-1)^T$. The matrices are given by
\begin{equation*}
    A_s = 
    \begin{pmatrix} 
        \lambda_s & 0 & 0 \\ 
        \lambda_s + a_s - b_s\sqrt{2} & -a_s & -b_s \\ 
        -(\lambda_s + a_s)\sqrt{2} - b_s & b_s & -a_s 
    \end{pmatrix}
    \qquad\text{and}\qquad
    A_u = 
    \begin{pmatrix} 
        a_u & -\lambda_u - a_u + b_u\sqrt{2} & -b_u \\ 
        0 & -\lambda_u & 0 \\ 
        b_u & -(\lambda_u + a_u)\sqrt{2} - b_u & a_u 
    \end{pmatrix},
\end{equation*}
with $\lambda_s, a_j > 0$ for $j\in\set{u,s}$, and $b_s b_u \neq 0$.

Now we proceed to prove that the requirements \ref{CB1}, \ref{CB2} and \ref{CB3} are satisfied.

\begin{enumerate}
	\item[\ref{CB1}] A direct calculation shows that the eigenvalues of $A_s$ are $\lambda_s$ and $-a_s \pm i b_s$, so $p_s$ is a hyperbolic saddle-focus for $X$. Similarly, the eigenvalues of $A_u$ are $-\lambda_u$ and $a_u \pm i b_u$, so $p_u$ is a hyperbolic saddle-focus for $Y$.
	
	\item[\ref{CB2}] We note that $W^s(p_s) = \set{(0,y,z)^T}$, so in a restricted neighborhood $U \subset \Sigma$, we have
$
W_\Sigma^s(p_s) = \set{(0, y, 0)^T \in U}.
$
Similarly, $W^u(p_u) = \set{(x,0,z)^T}$, thus,
$
W_\Sigma^u(p_u) = \set{(x, 0, 0)^T \in U},
$
giving us
$
p = (0,0,0)^T \in W_\Sigma^s(p_s) \pitchfork W_\Sigma^u(p_u).
$

To prove that $p\in\Sigma^c$, we simply note that
$
X(p) = (0, b_s, a_s)^T$ and $Y(p) = (-b_u, 0, a_u)^T,
$
therefore $\inner{X(p), \nabla f(p)} = a_s > 0$ and $\inner{Y(p), \nabla f(p)} = a_u > 0$, making both fields point into $R^+$.

Finally, the backward trajectory $\Gamma_{p_u \to p} = \set{\flow_Y^t(p) : t < 0} \subset R^-$, and the forward trajectory $\Gamma_{p \to p_s} = \set{\flow_X^t(p) : t > 0}\subset R^+$, as required.

	\item[\ref{CB3}] For the last condition, we note that the eigenvectors associated with the eigenvalues $\lambda_s$ and $-\lambda_u$ are $v_s = (1, 1, -\sqrt{2})^T$ and $v_u = (1, 1, \sqrt{2})^T$, respectively. The $1$-dimensional trajectories emerging from $p_s$ and $p_u$ along these directions can be parametrized as $p_s + t_s v_s$ and $p_u + t_u v_u$, and they intersect $\Sigma$ by setting $z=0$, giving us $t_s = 1/\sqrt{2}$ and $t_u = 1/\sqrt{2}$. This implies that both half-lines intersect $\Sigma$ at the exact same point
$
q = \left(\frac{1}{\sqrt{2}}, \; \frac{1}{\sqrt{2}}, \; 0\right)^T \in \Sigma.
$
Evaluating the vector fields at $q$ gives $\inner{X(q), \nabla f(q)} = -\lambda_s < 0$ and $\inner{Y(q),\nabla f(q)} = -\lambda_u < 0$. Both fields point into $R^-$, ensuring $q \in \Sigma^c$, $\Gamma_{p_s \to q} \subset R^+$, and $\Gamma_{q \to p_u} \subset R^-$.
\end{enumerate} 

Therefore the system in \eqref{eq:crossing_bykov} is a crossing Bykov cycle.
\end{example}

We emphasize that Example~\ref{example:crossing_bykov} falls under the hypotheses of Theorem~\ref{thm:crossing_bykov}, from which we can draw its conclusions.

\section{Sliding Bykov Cycles: Proof of Theorem \ref{thm:sliding_bykov}}\label{section:sliding_bykov_proof}

In this section we study a version of Bykov cycles intrinsic to Filippov systems, providing a structure to analyze it, in particular for proving Theorem~\ref{thm:sliding_bykov}.
From now on in this section, our system satisfies \ref{cond:linearization}--\ref{cond:distortion_bounds}.

In order to prove Theorem~\ref{thm:sliding_bykov}, we need to construct the first-return map that will be studied.

Our strategy will be to take a point in a small neighborhood of $q$ in $\partial\Sigma^s$ (the boundary of the sliding region), and follow it in the path induced by the flow of the Filippov system. This will allow us to separate its dynamics into local and global behaviours.

Let us consider, as given in Definition~\ref{def:sliding_bykov}, a sliding Bykov cycle $\Gamma := \set{p_u}\cup\Gamma_{p_u\to p_s}\cup\set{p_s}\cup\Gamma_{p_s\to p_u}$ with visible fold-regular point $q\in\Gamma_{p_u\to p_s}\cap\partial\Sigma^s$, and as in \eqref{eq:gamma}, for a sufficiently small $r>0$, the set
\[
\gamma = \gamma_r := \cl{B_r(q)\cap\partial\Sigma^s}
\]
is a compact curve containing $q$. The point $q$ separates this curve into two segments, and we can only guarantee that one of these segments follows closely the sliding Bykov cycle and meets $\Sigma$ again near $p_u$. By an abuse of notation, $\gamma$ will denote this specific segment from now on (including $q$ to ensure it remains a compact set).

Notice that we can apply a diffeomorphism $\gamma \to [0, r]$ (for instance, an arc-length parametrization), dealing with both interchangeably; all of the results in Theorem~\ref{thm:sliding_bykov} are invariant under this kind of parametrization, so we can treat the initial point as a local parameter $w \ge 0$, where the limit $w \to 0$ corresponds to $w \to q$.

We divide the proof of Theorem~\ref{thm:sliding_bykov} into the following steps.

\begin{enumerate}
    \item[\textit{Step 1.}]
    \textit{The Map from $\gamma$ to $\Sigma^s$.}
    We analyze the trajectories induced by the flow of $X$ starting at $\gamma$; in particular, we study the spiraling effects of the flow near $p_s$ using its linearization, and the resulting projected curve $\mu$ on $\Sigma^s$. 

    \item[\textit{Step 2.}]
    \textit{The Map from $\Sigma^s$ to $\gamma$.}
    Here, we consider the linearization of $Z_\Sigma$ and the family of trajectories $\Gamma_{[0,r]}$ in the past of $\gamma$ by the flow of $Z_\Sigma$.
    
    \item[\textit{Step 3.}]
    \textit{The Domain of the First-Return Map.}
    Then, we describe the intersection $\mu\cap\Gamma_{[0,r]}$ between the projected curve and the past flow of $\gamma$, from which we construct the domain of the first-return.
    
    \item[\textit{Step 4.}]
    \textit{The First Return Map as a CIFS.}
    We associate the first-return with a \emph{conformal iterated function system (CIFS)}, an analytical tool that will allow us to deduce several geometric properties of the invariant set.
    
    \item[\textit{Step 5.}]
    \textit{Proof of the Items of Theorem \ref{thm:sliding_bykov}.}
    Finally, we apply the theory of CIFS to deduce each one of the main results.
\end{enumerate}

\begin{figure}[h!]
\centering
\begin{overpic}[width=0.4\linewidth]{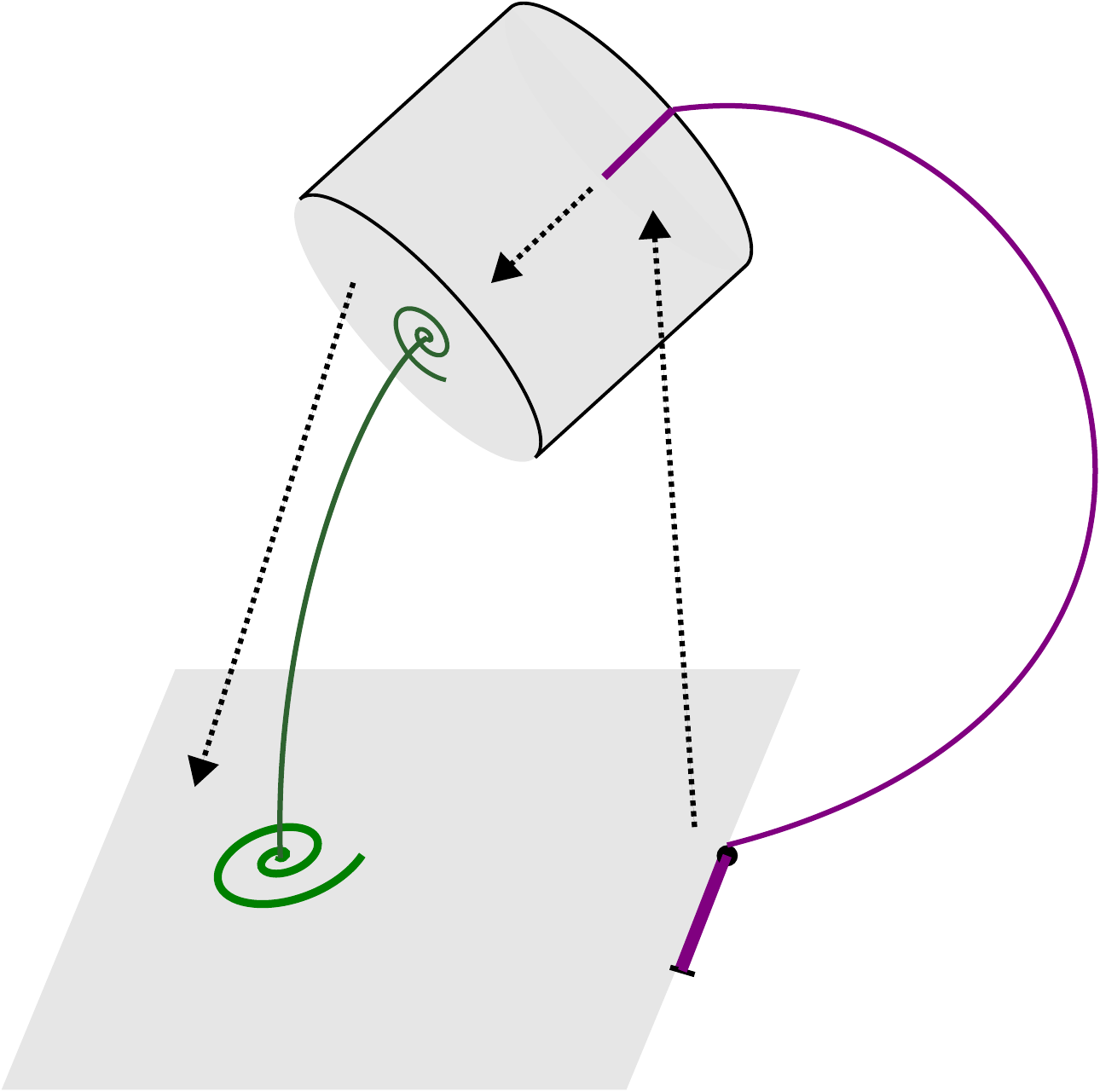}
\put(63,50){$\mu_{in}$}
\put(40,80){$\mu_{out}$}
\put(13,50){$\mu_{proj}$}
\put(67,15){$\gamma$}
\put(10,5){$\Sigma$}
\end{overpic}
\caption{Graphical representation of the intermediary maps of a sliding Bykov cycle.}
\label{fig:maps}
\end{figure}

So now we proceed with the proof.

\subsection{The Map from $\gamma$ to $\Sigma^s$}

The aim of this subsection is to deduce the following lemma. Here, an \emph{asymptotically periodic function $f$} is a function that can be written as $f=g+h$, where $g$ is periodic and $h(x)\to 0$ as $x\to\infty$.

\begin{lemma}\label{lemma:from_gamma_to_Sigma}
For a sufficiently small $r>0$, the image of $\gamma$ by the map induced by the flow of $X$ into $\Sigma^s\subset\Sigma$ is given by
\begin{equation}\label{eq:projected_curve_polar}
\hat{\mu}(t) : \left\{ \begin{aligned}
r_{\hat{\mu}}(t) &= \rho_{\hat{\mu}}(t) e^{-a_s t}, \\
\theta_{\hat{\mu}}(t) &= S b_s t + \omega_{\hat{\mu}}(t),
\end{aligned} \right.
\end{equation}
where $S = \sgn(\det(L)) \in \{-1, 1\}$, with $L$ given by \eqref{eq:L_definition}, and the amplitude $\rho_{\hat{\mu}}(t) > 0$ and phase shift $\omega_{\hat{\mu}}(t)$ are bounded, asymptotically periodic functions.
\end{lemma}

To prove it, we will analyze the intermediate maps from $\gamma$ to a neighborhood of $p_s$, the local map near $p_s$, and the projection from the neighborhood of $p_s$ to $\Sigma^s$. We remark that when analyzing the local map near $p_s$ and the resulting projected curve into $\Sigma^s$ (near $p_u$), we are going to use linearized coordinates.

So first, by \ref{cond:linearization} there exist a neighborhood $\mathcal{U}(p_s)$ near $p_s$ such that the flow is $C^2$-conjugate to a linear flow by the map $h_s: \mathcal{U}(p_s) \to \mathbb{R}^3$.

We then take $V^{s}\subset\mathcal{U}(p_s)$ a closed neighborhood for $p_s$ such that
\[
\hat{V}^{s} := h_s(V^{s}) = \set{(r,\theta,z):r\in(0,r_{in}], z\in[0,z_{out}]} \cup \set{(r,z):r=0, z\in[0,z_{out}]},
\]
is a cylinder with cylinder wall and top disk given by
\[
\hat{V}_{in}^{s} := \set{(r, \theta, z) : r = r_{in}, z \in [0, z_{out}]}
\]
and
\[
\hat{V}_{out}^{s} = \set{(r, \theta, z) : r \in (0, r_{in}], z = z_{out}} \cup \set{(r,z):r=0, z=z_{out}}.
\]

Now we consider an initial condition $w \in \gamma$ near $q$. We recall that we interchange $w$ between the point on the curve $\gamma$ and the arc-length parametrization on $[0, r]$.

The trajectory under the flow of $X$ starts at $w \in \gamma \subset \partial\Sigma^s$ near the fold point $q$ (where $w = 0$ corresponds to $q$). The trajectory crosses the cylinder wall $V_{in}^{s}$ at a flight time $t_{in}(w) > 0$.

We then define the $C^2$ map $\mu_{in}:\gamma\to V_{in}^{s}$ as
\[
\mu_{in}(w) := \flow_X^{t_{in}(w)}(w),
\]
where $t_{in} := \inf\set{t>0 : \flow_X^t(w) \in V_{in}^{s}}$, and we denote $\mu_{in} := \mu_{in}(\gamma) \subset V_{in}^{s}$.

The linearized version of this map is the $C^2$ map $\hat{\mu}_{in}:\gamma\to\hat{V}_{in}$ given by
\begin{equation*}
\hat{\mu}_{in}(w) := h_s(\mu_{in}(w)),
\end{equation*}
and we denote $\hat{\mu}_{in} := h_s(\mu_{in}) = \hat{\mu}_{in}(\gamma) \subset \hat{V}_{in}^{s}$.

The cylindrical coordinates of the points in the curve $\hat{\mu}_{in}$ are thus given by
\begin{equation}\label{eq:entry_r}
\hat{\mu}_{in}(w) : 
	\begin{cases}
	r_{in}(w) = r_{in},										\\
	\theta_{in}(w) = \theta_0 + c_{in} w + O(w^2),	\\
	z_{in}(w) = d_{in} w + O(w^2),
	\end{cases}
\end{equation}
where $d_{in} > 0$ and $\theta_0$ is the angle associated with the point $q$.

Inside the neighborhood $\hat{V}^{s}$, the local flow of $X$ is conjugate via $h_s$ to the linear system
\begin{equation*}
\dot{r} = -a_s r,	\qquad
\dot{\theta} = b_s,	\qquad
\dot{z} = \lambda_s z.
\end{equation*}
which has the solution for the initial value in \eqref{eq:entry_r} given by
\begin{equation*}
r(t,w) = r_{in} e^{-a_s t},	\qquad
\theta(t,w) = b_s t + \theta_{in}(w),	\qquad
z(t,w) = z_{in}(w)e^{\lambda_s t}.
\end{equation*}

The solution for the local flight time $t(w)$ required to travel from the cylinder wall $\hat{V}_{in}^{s}$ to the exit disk $\hat{V}_{out}^{s}$ depends only on $z_{in}(w)$ and is given by
\begin{align}\label{eq:t_w}
t(w) 
= \frac{1}{\lambda_s} \log\left( \frac{z_{out}}{d_{in}w(1+O(w))} \right) &= t_0 -\frac{1}{\lambda_s} \log(w) - \frac{1}{\lambda_s}\log(1+ O(w))
= t_0 -\frac{1}{\lambda_s} \log(w) +  O(w) 
\end{align}
with $t_0 = \frac{1}{\lambda_s}\log\left(\frac{z_{out}}{d_{in}}\right)$. Notice that
\begin{equation}\label{eq:derivative_t_w}
\frac{dt}{dw}(w) =
-\frac{1}{\lambda_s w} + O(1) \implies 
\frac{c_1}{w} \leq \abs{\frac{dt}{dw}(w)} \leq \frac{c_2}{w},
\end{equation}
for some $c_1,c_2>0$.

For sufficiently small $w>0$ we have $dt/dw < 0$, allowing us to invert the relationship in \eqref{eq:t_w}, and since as $w\to0$ we have $e^{O(w)} = 1+ O(w)$, we obtain
\begin{equation}\label{eq:w(t)}
w = w(t) = K_w e^{-\lambda_s t}(1 + O(w(t))) \implies w(t) = K_w e^{-\lambda_s t}+O(e^{-2\lambda_s t}) = K(t)e^{-\lambda_s t},
\end{equation}
where $K_w := e^{\lambda_s t_0}$ and $K(t)$ is a bounded function. Besides, by inverting \eqref{eq:derivative_t_w}, we also have
\begin{equation}\label{eq:derivative_w(t)}
\frac{dw}{dt}(t) = -\lambda_s K_w e^{-\lambda_s t} + O(e^{-2\lambda_s t}) \implies
C_1w\leq \abs{\frac{dw}{dt}(t)} \leq C_2w,
\end{equation}
for some $C_1,C_2>0$.

From this, we can define the map $\hat{\mu}_{out}:\hat{\mu}_{in}\setminus\set{\hat{\mu}_{in}(q)}\subset \hat{V}_{in}^{s} \to \hat{V}_{out}^{s}$ given by
\[
\hat{\mu}_{out}(\hat{\mu}_{in}(w)) = \hat{\mu}_{out}(\hat{\mu}_{in}(w(t))) := \begin{pmatrix}
				r_{in} e^{-a_s t} \\
				b_s t + \theta_{in}(w(t)) \\
				z_{out}
				\end{pmatrix}.
\]

Again, we denote the curve $\hat{\mu}_{out} := \hat{\mu}_{out}(\hat{\mu}_{in})$. As we can see in the previous formula, we can parametrize it by the time $t$; in Cartesian coordinates it is written as
\begin{align*}\label{eq:hat_mu_out}
\hat{\mu}_{out}(t) = \begin{pmatrix}
				0 \\
				0 \\
				z_{out}
				\end{pmatrix}
				+ r_{in} e^{-a_s t}\left[\begin{pmatrix}
				\cos(b_s t + \theta_0) \\
				\sin(b_s t + \theta_0) \\
				0
				\end{pmatrix}
				+ O(e^{-\lambda_s t})\right] 
				= \hat{u}_{out}
				+ r_{in} e^{-a_s t}[\hat{v}_0(t)
				+ O(e^{-\lambda_s t})],
\end{align*}
where $\hat{u}_{out} := (0, 0, z_{out})^T$ and $\hat{v}_0(t) := (\cos(b_s t + \theta_0), \sin(b_s t + \theta_0), 0)^T$.

Notice that
\begin{equation*}\label{eq:derivative_mu_out}
\norm{\frac{d}{dt}\hat{\mu}_{out}(t)} = O(e^{-a_st}) \implies \norm{\frac{d}{dw}\hat{\mu}_{out}(w)} = O(e^{-a_st}) \abs{\frac{dt}{dw}}.
\end{equation*}

In order to go back to the original coordinate system, we define
\[
\mu_{out}(t) = h_s^{-1}(\hat{\mu}_{out}(t)).
\]

As $t \to \infty$ (corresponding to $w \to 0$), the point $\hat{\mu}_{out}(t)$ approaches the center of the top disk , which lies on the linearized local unstable manifold. Let $u_{out} := h_s^{-1}(\hat{u}_{out})$ be the corresponding point on the unstable heteroclinic connection $\Gamma_{p_s \to p_u}.$

Now we describe the transition map $\mu_{proj} : V_{out}^{s} \to \Sigma$, which takes a neighborhood of $u_{out}$ and projects it along the trajectories of $X$ onto a neighborhood of the pseudo-saddle-focus $p_u \in \Sigma^s$. More precisely, it is defined as
\[
\mu_{proj}(w) := \flow_X^{t_{proj}(w)}(w), 
\]
where $t_{proj} := \inf\set{t>0 : \flow_X^t(w)\in\Sigma}$. Notice that since $u_{out}\in\Gamma_{p_s \to p_u}$, then $\mu_{proj}(u_{out}) = p_u$.

Finally, we define the linearized global transition map $\hat{\mu}_{glob}: \hat{V}^{s}_{out} \to \mathbb{R}^2$ directly between the two linearized coordinate spaces as the composed diffeomorphism
\begin{equation*}
\hat{\mu}_{glob} := h_u \circ \mu_{proj} \circ h_s^{-1}.
\end{equation*}
Since each individual component is a local $C^2$ diffeomorphism, $\hat{\mu}_{glob}$ is a local $C^2$ diffeomorphism.

We then consider the linear map $L: \mathbb{R}^2 \to \mathbb{R}^2$ given by
\begin{equation}\label{eq:L_definition}
L := D\hat{\mu}_{glob}(\hat{u}_{out}) = Dh_u(p_u) \circ D\mu_{proj}(u_{out}) \circ D(h_s^{-1})(\hat{u}_{out}).
\end{equation}

We can now write the projected curve in the linearized coordinates of $p_u$, denoted by $\hat{\mu}(t) := \hat{\mu}_{glob}(\hat{\mu}_{out}(t))$. We can write $\hat{\mu}$ around $\hat{u}_{out}$ as
\begin{equation}\label{eq:projected_curve_with_error}
\hat{\mu}(t) = \hat{\mu}_{glob}(\hat{u}_{out}) + D\hat{\mu}_{glob}(\hat{u}_{out})(\hat{\mu}_{out}(t) - \hat{u}_{out}) + O(\norm{\hat{\mu}_{out}(t) - \hat{u}_{out}}^2) = r_{in} e^{-a_s t}L[\hat{v}_0(t) + O(e^{-\delta t})],\end{equation}
where $\delta := \min\set{a_s, \lambda_s}$.

Near the pseudo-saddle-focus we can represent $\hat{\mu}(t)$ in local polar coordinates as
\begin{equation*}
r_{\hat{\mu}}(t) = \rho_{\hat{\mu}}(t) e^{-a_s t}, \qquad
\theta_{\hat{\mu}}(t) = S b_s t + \omega_{\hat{\mu}}(t),
\end{equation*}
completing the proof of Lemma \ref{lemma:from_gamma_to_Sigma}.

\subsection{The Map from $\Sigma^s$ to $\gamma$}

In this section we prove the following lemma.

\begin{lemma}\label{lemma:from_Sigma_to_gamma}
Consider the trajectories induced by $Z_\Sigma$, denoted $\Gamma_0(\tau)$, which reaches $q \in \gamma$, and $\Gamma_r(\tau)$, which reaches $w=r$, the opposite extremity of $\gamma$. 
Then for a neighborhood of $p_u$, we can write them in linearizing coordinates as
\begin{equation}\label{eq:B_nu}
\hat{\Gamma}_j(\tau):
\begin{cases}
r_{\hat{\Gamma}_j}(\tau) = \eta e^{-a_u \tau},	\\
\theta_{\hat{\Gamma}_j}(\tau) = \hat{\Theta}_j - b_u \tau,
\end{cases}
\end{equation}
for some constants $\eta$ and $\hat{\Theta}_j$, with $j \in \set{0,r}$.
\end{lemma}

Similarly to what was done for $p_s$, by \ref{cond:linearization} there exist a neighborhood $\mathcal{U}(p_u)$ near $p_u$ such that the flow is $C^2$-conjugate to a linear flow by the map $h_u: \mathcal{U}(p_u) \to \mathbb{R}^2$.

We then take $V_{\eta}^{u} \subset \mathcal{U}(p_u) \subset \Sigma^s$ a closed neighborhood for $p_u$ defined by
\begin{equation*}
\hat{V}_{\eta}^{u} := h_u(V_{\eta}^{u}) = \cl{B_{\eta}(p_u)} = \set{(r,\theta):r\in(0,\eta]} \cup \set{0}.
\end{equation*}

In this coordinate system, the linearized sliding vector field $\hat{Z}_{\Sigma}$ is written as
    \begin{equation*}
    \dot{w} = \hat{Z}_{\Sigma}(w) = \hat{A}_u w, \qquad \text{where } \hat{Z}_{\Sigma}(w_1, w_2) = (a_u w_1 - b_u w_2, b_u w_1 + a_u w_2) \iff \hat{A}_u = \begin{pmatrix}
    	a_u & -b_u \\
    	b_u & a_u
    \end{pmatrix}.
    \end{equation*}

Since we are interested in the set of points that eventually reach $\gamma$ under the flow of $Z_{\Sigma}$, it is natural to consider the set
\[
\Gamma_{[0,r]} := \bigcup_{t\le0}\flow_{Z_{\Sigma}}^t(\gamma) = \bigcup_{s\in[0,r]}\Gamma_s,\qquad\text{where }\Gamma_s := \bigcup_{t\le0}\set{\flow_{Z_{\Sigma}}^t(s)}.
\]
consisting of points of $\Sigma$ that reach $\gamma$ after some positive time.

We now take the linearized spirals $\hat{\Gamma}_{j} := h_u(\Gamma_{j}) \cap \hat{V}_{\eta}^{u}$, for $j\in\set{0,r}$, and define $\hat{\Theta}_0$ and $\hat{\Theta}_r$, respectively, as the coordinate angles of the endpoints $\hat{\Gamma}_{0}\cap\partial \hat{V}_{\eta}^{u}$ and $\hat{\Gamma}_{r}\cap\partial \hat{V}_{\eta}^{u}$.

Notice that we could parametrize $\hat{\Gamma}_{j}$ by the positive sliding time $\tau$ required to reach $\partial V_{\eta}^{u}$, getting
\begin{equation*}
r_{\hat{\Gamma}_j}(\tau) = \eta e^{-a_u \tau},	\qquad
\theta_{\hat{\Gamma}_j}(\tau) = \hat{\Theta}_j - b_u \tau,
\end{equation*}
for $j \in \set{0, r}$, just as stated in Lemma \ref{lemma:from_Sigma_to_gamma}.

\subsection{The Domain of the First-Return Map}

In this section, we are going to prove the next lemma.
\begin{lemma}\label{ppst:domain_first_return}
The set of points in $\gamma$ with a well-defined first-return is given by the disjoint union of closed intervals
\[
W^{\gamma} = W_r^{\gamma} := \bigcup_{k\ge1} I_k^{\gamma}\subset\gamma,
\]
with $I_k^\gamma:=[w_k^{\min},w_k^{\max}]$.
Furthermore, for $w_k\in I_k^\gamma$, there exist constants $C_{min}, C_{max} > 0$ such that
$
C_{min} \exp\left( - \frac{2\pi k}{|\kappa|} \right) \le w_k \le C_{max} \exp\left( - \frac{2\pi k}{|\kappa|} \right).
$
\end{lemma}

\begin{figure}[h!]
\centering
\vspace{0.2cm}
\begin{overpic}[width=0.8\linewidth]{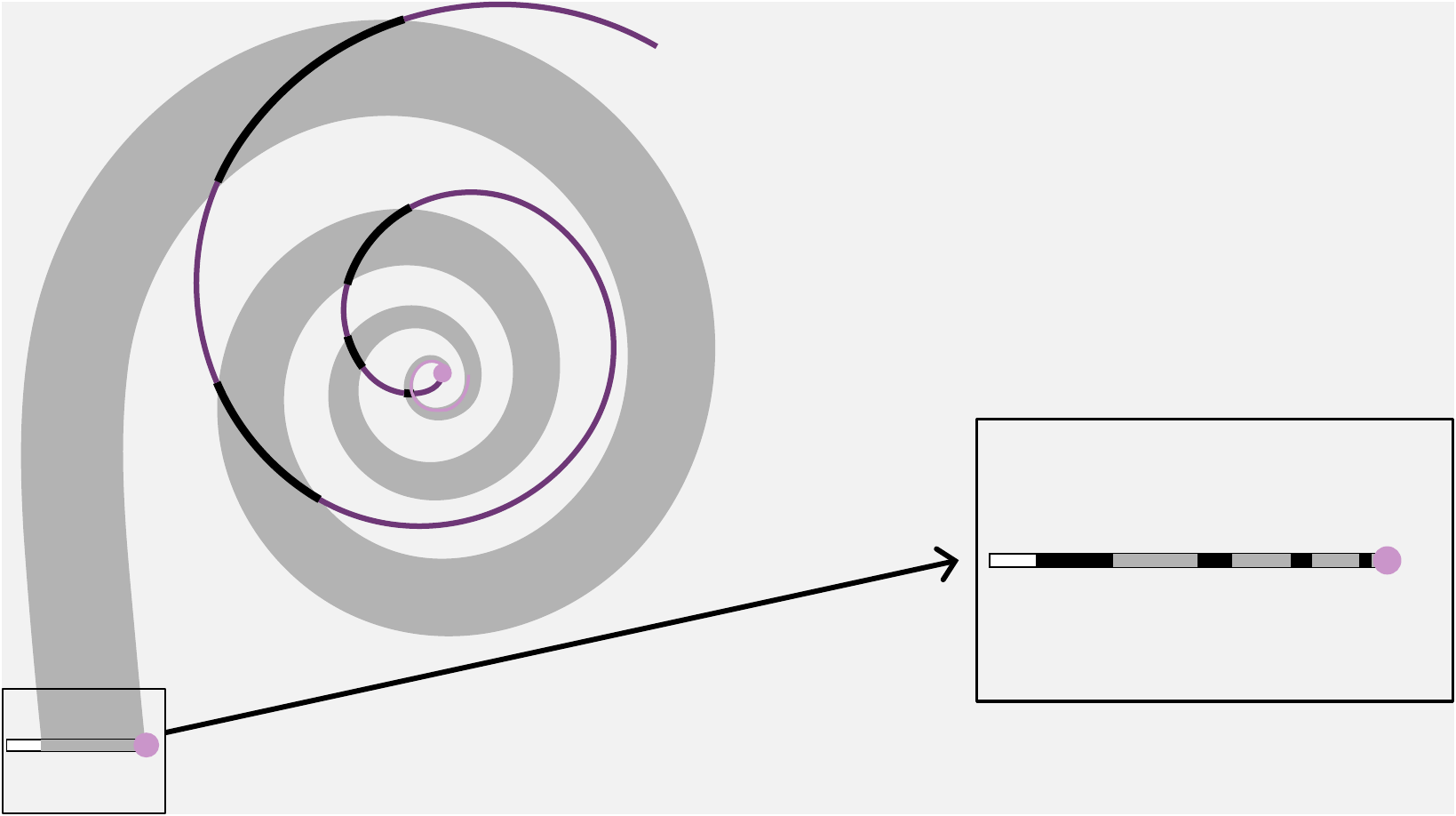}
\put(5,2){$\gamma_r$}
\put(10,2){$q$}

\put(46,50){$\mu_r$}

\put(16,50){$I_1^{\mu_r}$}
\put(18,19){$I_2^{\mu_r}$}
\put(22,39){$I_3^{\mu_r}$}
\put(22,28){$I_4^{\mu_r}$}

\put(31.5,30){$p_u$}

\put(2.5,20){$\Gamma_{[0,r]}$}

\put(82,12){$\gamma_r$}

\put(73,19){$I_1^{\gamma_r}$}
\put(83,19){$I_2^{\gamma_r}$}
\put(89,19){$I_3^{\gamma_r}$}
\put(93,19){$I_4^{\gamma_r}$}

\put(95,14){$q$}

\end{overpic}
\caption{Representation of the construction of $W^{\mu}$ and $W^{\gamma}$.}
\label{fig:intervals}
\end{figure}

Since $\mu\setminus\set{p_u}$ near $p_u$ is transverse to $Z_{\Sigma}$ (see Proposition \ref{pps:F_greater_zero} from Appendix~\ref{appendix:transversality}) and $p_u\in\mu\subset \Sigma^s$ is a hyperbolic focus of $Z_{\Sigma}$, we deduce that for $r>0$ sufficiently small we have
\[
W^{\mu} = W_r^{\mu} := \mu\cap \Gamma_{[0,r]} = \bigcup_{i\ge1} I_i^{\mu} \subset \mu
\]
where $\set{I_i^{\mu}}_{i\ge1}$ are the one-dimensional connected components of $W^{\mu}$. They satisfy $I_i^{\mu}\cap I_j^{\mu} = \emptyset$ for $i\neq j$, and $I_i^{\mu} \to \set{p_u}$ in the Hausdorff distance as $i\to\infty$ (see Figure~\ref{fig:intervals}). We define the index set
\[
\cal{J}^{\mu} = \cal{J}^{\mu_r} := \set{I_i^{\mu}}_{i\ge1}.
\]

Similarly, for each $i\ge1$, we define
\[
I_i^{\gamma} := \theta_X^{-1}(I_i^{\mu}) \subset \gamma
\]
where $\theta_X:\gamma\setminus\set{q}\to\mu\setminus\set{p_u}$ is the diffeomorphism induced by the flow of $X$ (again, see Figure~\ref{fig:intervals}). This ensures that the collection $\set{I_i^\gamma}_{i\ge1}$ has the same properties $I_i^{\gamma}\cap I_j^{\gamma} = \emptyset$, if $i\neq j$. Besides, $I_i^{\gamma} \to \set{q}$ in the Hausdorff distance as $i\to\infty$. Thus, we define
\[
W^{\gamma} = W_r^{\gamma} := \bigcup_{i\ge1} I_i^{\gamma}\subset\gamma
\]
and the index set
\[
\cal{J} = \cal{J}^{\gamma} = \cal{J}^{\gamma_r} := \set{I_i^{\gamma}}_{i\ge1}.
\]

Using Lemma \ref{lemma:from_gamma_to_Sigma} and Lemma \ref{lemma:from_Sigma_to_gamma}, we have that the curve $\hat{\mu} = \hat{\mu}(t)$ and the boundary $\hat{\Gamma}_j = \hat{\Gamma}_j(\tau)$ intersect if there exist a flight time $t$, a sliding time $\tau$, and an integer $k \in \mathbb{Z}$ such that their coordinates match, that is, if the following conditions are simultaneously met,
\begin{equation}\label{eq:equality}
\begin{cases}
r_{\hat{\mu}}(t) = r_{\hat{\Gamma}_j}(\tau),		\\
\theta_{\hat{\mu}}(t) = \theta_{\hat{\Gamma}_j}(\tau) + 2\pi k.
\end{cases}
\end{equation}

Substituting the radial parameterizations in \eqref{eq:projected_curve_polar} and \eqref{eq:B_nu} in the radial equality condition~\eqref{eq:equality} we get $\rho_{\hat{\mu}}(t) e^{-a_s t} = \eta e^{-a_u \tau}$, from which we deduce
\begin{equation}\label{eq:tau}
\tau = \tau(t) = \frac{a_s}{a_u} t - \frac{1}{a_u} \log(\rho_{\hat{\mu}}(t)) + \frac{1}{a_u} \log(\eta).
\end{equation}

This allows us to define the map $\hat\nu_{\eta}:\hat{\mu}\to \partial \hat V_{\eta}^u$ given by
\begin{equation*}\label{eq:pi_nu}
\hat\nu_{\eta}(\hat{\mu}(t)) := \left(\flow_{\hat{Z}_{\Sigma}}^{\tau(t)}(\hat{\mu}(t))\right)
\end{equation*}
the function that maps a point in $\hat{\mu}$ to its image under the flow of $\hat{Z}_{\Sigma}$ in the boundary of $\hat V_{\eta}^u$. In order to view it in non-linearized coordinates, we define $\nu_{\eta}:\hat{\mu}\to \partial V_{\eta}^u$ as $\nu_{\eta} = h_u^{-1}\circ\hat\nu_{\eta}$.

Moreover, we now can also define $\nu_{Z_{\Sigma}}:\Gamma_{[0,r]}\cap\partial V_{\eta}^u\to\gamma$ as 
\begin{equation*}
\nu_{Z_{\Sigma}}(x) := \flow_{Z_{\Sigma}}^{t_{Z_{\Sigma}}(x)}(x),
\end{equation*}
where
$t_{Z_{\Sigma}}(x) := \inf\set{t>0 : \flow_{Z_{\Sigma}}^t(x) \in\gamma}$.

Substituting the expression for $\tau$ in \eqref{eq:tau} into the angular equality condition~\eqref{eq:equality} gives us
\begin{equation*}
S b_s t + \omega_{\hat{\mu}}(t) = \hat{\Theta}_j - b_u \left( \frac{a_s}{a_u} t - \frac{1}{a_u} \log(\rho_{\hat{\mu}}(t)) + \frac{1}{a_u} \log(\eta) \right) + 2\pi k,
\end{equation*}
and by putting together all linear $t$ terms and absorbing bounded terms into $B(t) = O(1)$, we have
\begin{equation*}
\left( S b_s + b_u \frac{a_s}{a_u} \right) t = 2\pi k + \hat{\Theta}_j + B(t)
\end{equation*}

We now define $\kappa := \frac{1}{\lambda_s}\left( S b_s + b_u \frac{a_s}{a_u} \right)$, which by  \ref{cond:non_resonance_undistorted} we have $\kappa\neq 0$ (by multiplying $a_u/b_ub_s \neq 0$ while recalling that $c_s = -a_s/b_s$). By choosing the indexing of the intersection times according to the
orientation of the orbit, we may assume that $t_k^{(j)}\to+\infty$ as
$k\to+\infty$. Thus, after possibly replacing $k$ by $-k$ when
$\kappa<0$, the intersection times satisfy
\begin{equation}\label{eq:intersection_times}
t_k^{(j)} = \frac{2\pi k+\hat\Theta_j+B(t_k^{(j)})}
{\abs{\kappa}\lambda_s}.
\end{equation}

Applying the intersection times given by \eqref{eq:intersection_times} in \eqref{eq:w(t)} gives the endpoints
\[
w_k^{(0)} := w(t_k^{(0)}) \qquad\text{ and }\qquad w_k^{(r)} := w(t_k^{(r)})
\]
of the interval $I_k^\gamma:=[w_k^{\min},w_k^{\max}]$, where we are denoting
$
w_k^{\min} := \min\set{w_k^{(0)},w_k^{(r)}}$
and
$
w_k^{\max} := \max\set{w_k^{(0)},w_k^{(r)}},
$
in such a way that, for $w_k\in I_k^\gamma$, we deduce that there exist constants $C_{min}, C_{max} > 0$ with
\begin{equation}\label{eq:bounds_w_j_k}
C_{min} \exp\left( - \frac{2\pi k}{|\kappa|} \right) \le w_k \le C_{max} \exp\left( - \frac{2\pi k}{|\kappa|} \right),
\end{equation}
as desired.

Notice that, as $r\to 0$, we have $I_1^{\gamma} = I_1^{\gamma_r}\to\set{q}$.

\subsection{The First Return Map as a CIFS}

We define the first-return
\[
\pi(w)=\pi|_{W^{\gamma}}(w)
:=\nu_{Z_\Sigma}
\circ h_u^{-1}
\circ \hat\nu_{\eta}
\circ \hat{\mu}_{glob}
\circ \hat\mu_{out}
\circ \hat{\mu}_{in}
\big|_{W^{\gamma}}(w):W^{\gamma}\rightarrow\gamma
\]
and for $J_k:=I_k^\gamma\subset W^\gamma$ we define
\begin{equation}\label{eq:psi}
\psi_{J_k}:=\left(\pi|_{J_k}\right)^{-1} : \gamma\to J_k
\qquad\text{ and }\qquad
\Psi := \set{\psi_{J_k}}_{k\in\N}.
\end{equation}

Our purpose in this section is to prove that the set $\Psi$ is a \emph{conformal iterated function system (CIFS)}, which will allow us to deduce the equality between the invariant set given in \eqref{eq:local} and the attractor set in \eqref{eq:atractor_set}, in order to deduce the results of Theorem~\ref{thm:sliding_bykov}. We give an introductory presentation of CIFS theory in Appendix~\ref{appendix:cifs}.

\subsubsection{$\Psi$ as a CIFS}
We now prove the next lemma. We recall that the conditions for a IFS to be conformal are given in Appendix~\ref{appendix:cifs}.

\begin{lemma}\label{lemma:cifs}
For a sufficiently small $r>0$, the family of maps $\Psi$ in \eqref{eq:psi} is a conformal iterated function system.

\end{lemma}

\begin{proof}
The most immediate conditions are \ref{conf:c1} and \ref{conf:c3}, which are direct from the construction of $W^\gamma$ and $\psi_{J_k}$; then condition \ref{conf:c5} is deduced from the fact that $\gamma = [0,r]$ is a closed interval, so the expression is equal to $r/2$. Also, provided that each $\psi_{J_k}$ is a contraction, which will be proved soon for \ref{conf:c2}, the condition \ref{conf:c4} is valid since, for some $\delta>0$, we could consider the interval $I = I^\delta := (-\delta, r+\delta)\supset [0,r] = \gamma$ and simply extend in a $C^2$ fashion each map $\psi_{J_k}:\gamma\to\gamma$ to $\psi_{J_k}^I:I\to I$, which can be done in the one-dimensional context.

Now we focus on proving \ref{conf:c2}.
For that we consider the $t$-parametrization of $\pi$ given by
\[
\Pi(t) := \pi(w(t)) = \nu_{Z_{\Sigma}} \circ h_u^{-1} \circ \hat\nu_{\eta}
\circ \hat{\mu}(t),
\]
which gives the relation
\[
\frac{d}{dw}\pi(w) = \Pi'(t) \frac{dt}{dw}.
\]

Since $\nu_{Z_\Sigma}$ and $h_u^{-1}$ are $C^2$-diffeomorphisms on compact sets, there exist constants
$K_1, K_2>0$, such that
\[
K_1\cdot
\norm{
\frac{d}{dt}\hat\nu_{\eta}(\hat{\mu}(t))
}\cdot
\abs{\frac{dt}{dw}}
\leq
|\pi'(w)|
\leq
K_2\cdot
\norm{
\frac{d}{dt}\hat\nu_{\eta}(\hat{\mu}(t))
}\cdot
\abs{\frac{dt}{dw}}.
\]

So now we focus on the term $\norm{\frac{d}{dt}\hat\nu_{\eta}(\hat{\mu}(t))}$. By using
$
\hat\nu_{\eta}(r,\theta)
=
\left(
\eta,\,
\theta+\frac{b_u}{a_u}\log\frac{\eta}{r}
\right)
$
and applying it in $\hat{\mu}(t)$, which has radial and angular coordinates given by
$
r_{\hat\mu}(t)=\rho_{\hat\mu}(t)e^{-a_st}
$
and
$
\theta_{\hat\mu}(t)
=
Sb_st+\omega_{\hat\mu}(t),
$
we obtain
\begin{equation}\label{eq:nu_eta}
\hat\nu_{\eta}(\hat\mu(t))
=
\left(
\eta,\,
Sb_st+\omega_{\hat\mu}(t)
-\frac{b_u}{a_u}\log\rho_{\hat\mu}(t)
+\frac{b_ua_s}{a_u}t
+\frac{b_u}{a_u}\log\eta
\right).	
\end{equation}

From this, we have
\[
\norm{\frac{d}{dt}\hat\nu_\eta(\hat\mu(t))}
=
\norm{\left(0,G(t)\right)}
=
\abs{G(t)},
\]
where
$
G(t)
=
Sb_s+\frac{a_sb_u}{a_u}
+\omega'_{\hat\mu}(t)
-\frac{b_u}{a_u}
\frac{\rho'_{\hat\mu}(t)}{\rho_{\hat\mu}(t)}
=
\theta_{\hat\mu}'(t)
-\frac{b_u}{a_u}
\frac{r_{\hat\mu}'(t)}{r_{\hat\mu}(t)}.
$

On the other hand, we have the quantity
\[
F(t)
=
\det\left[
\hat\mu'(t),\hat Z_\Sigma(\hat\mu(t))
\right]\\
=
b_ur_{\hat\mu}(t)r_{\hat\mu}'(t)
-a_ur_{\hat\mu}(t)^2\theta_{\hat\mu}'(t)
=
-a_ur_{\hat\mu}(t)^2
\left(
\theta_{\hat\mu}'(t)
-\frac{b_u}{a_u}
\frac{r_{\hat\mu}'(t)}{r_{\hat\mu}(t)}
\right)
\]
(this quantity appears in more detail in Appendix~\ref{appendix:transversality}, in particular \eqref{eq:transversality_determinant}).
Thus,
\[
F(t)=-a_ur_{\hat\mu}(t)^2G(t).
\]

By Proposition~\ref{pps:F_greater_zero}, there exist
\(c_F,C_F>0\) and \(T>0\) such that
\[
c_F e^{-2a_st}
\le |F(t)|
\le C_F e^{-2a_st},
\]
for all $t\ge T$. Moreover, since
\[
r_{\hat\mu}(t)
=
\rho_{\hat\mu}(t)e^{-a_st},
\]
and \(\rho_{\hat\mu}\) is bounded above and bounded away from zero for
\(t\ge T\), there exist \(0<c_\rho\le C_\rho\) such that
\[
c_\rho e^{-a_st}
\le r_{\hat\mu}(t)
\le C_\rho e^{-a_st}.
\]
Therefore,
\[
\frac{c_F}{a_u C_\rho^2}
\le
|G(t)|
\le
\frac{C_F}{a_u c_\rho^2},
\]
for $t\ge T$. Consequently, $0<c\leq\norm{\frac{d}{dt}\hat\nu_\eta(\hat\mu(t))}\leq C$, for some $c,C>0$, which implies that there exist constants $K_3,K_4>0$, independent of $w$,
such that
\begin{equation}\label{eq:bounds_derivative_pi}
K_3\cdot
\abs{\frac{dt}{dw}(w)}
\leq
|\pi'(w)|
\leq
K_4\cdot
\abs{\frac{dt}{dw}(w)}.
\end{equation}

It follows from~\eqref{eq:derivative_t_w} that there exist constants $K_5,K_6>0$,
independent of $w$, such that for every sufficiently small $w$, we have
\begin{equation}\label{eq:derivative_pi_bounded_by_dt_dw}
\frac{K_5}{w}
\leq
|\pi'(w)|
\leq
\frac{K_6}{w}.
\end{equation}

For $z=\pi(w)\in J_k$, we conclude by~\eqref{eq:derivative_pi_bounded_by_dt_dw} that
\begin{equation}\label{eq:derivative_psi_z}
\frac{w_k^{\min}}{K_6}
\leq
\frac{w}{K_6}
\leq
|\psi_{J_k}'(z)|
\leq
\frac{w}{K_5}
\leq
\frac{w_k^{\max}}{K_5},
\end{equation}
for $z\in J_k$. Finally, by~\eqref{eq:bounds_w_j_k} and~\eqref{eq:derivative_psi_z} we have that $\psi_{J_k}$ as in \eqref{eq:psi} satisfies
\begin{equation*}\label{eq:bound_on_derivative_psi}
\frac{C_{\min}}{K_6}
\exp\left(-\frac{2\pi k}{|\kappa|}\right)
\leq
|\psi_{J_k}'(z)|\leq
\frac{C_{\max}}{K_5}
\exp\left(-\frac{2\pi k}{|\kappa|}\right),
\end{equation*}
which shows that each $\psi_{J_k}$ is a diffeomorphism and $\Psi$ is a family of uniform contractions, if $r > 0$ is chosen sufficiently small.

Finally, we need to prove \ref{conf:c6}. For that we first establish a second-derivative estimate for
the first-return map. We use the flight time $t$ as a parameter and
write
\[
\Pi(t)
=
\nu_{Z_\Sigma}\circ h_u^{-1}\circ\hat\nu_{\eta}\circ\hat{\mu}_{glob}
\circ\hat\mu_{out}\circ\hat\mu_{in}(w(t)).
\]
The maps $\nu_{Z_\Sigma}$, $h_u^{-1}$, $\hat{\mu}_{glob}$, and
$\hat\mu_{in}$ are $C^2$ diffeomorphisms from their domains onto their images, which are compact sets, so their first and second derivatives are uniformly bounded. It remains to analyze $\hat\mu_{out}$ and $\hat\nu_{\eta}$.

For the local map $\hat\mu_{out}$, we use its
$t$-parametrization. By (9),
\[
\hat\mu_{out}(t)
=
\hat u_{out}
+
r_{in}e^{-a_st}
\left[\hat v_0(t)+O(e^{-\lambda_st})\right].
\]
Lemma~\ref{lem:taylor_remainder_curve} allows us to estimate the remainder and conclude that
\[
\hat\mu_{out}'(t)=O(1)
\qquad\text{and}\qquad
\hat\mu_{out}''(t)=O(1).
\]
Since $\hat\mu(t)=\hat{\mu}_{glob}
(\hat\mu_{out}(t))$, the same estimates hold for
$\hat\mu(t)$.

Now we consider the sliding map $\hat\nu_{\eta}$. By
\eqref{eq:nu_eta}, its image is parametrized by
\[
r_{\hat\nu_{\eta}}(t)=\eta
\qquad\text{ and }\qquad
\theta_{\hat\nu_{\eta}}(t)
=
\left(Sb_s+\frac{b_ua_s}{a_u}\right)t
+\omega_{\hat\mu}(t)
-\frac{b_u}{a_u}\log\rho_{\hat\mu}(t)
-\frac{b_u}{a_u}\log\eta.
\]
Since $\rho_{\hat\mu}$ is bounded away from zero and
$\rho_{\hat\mu}$, $\omega_{\hat\mu}$ have bounded first and
second derivatives, the first and second derivatives of
$\theta_{\hat\nu_{\eta}}(t)$ are bounded. Hence the first and second
derivatives of the $t$-parametrized image under
$\hat\nu_{\eta}$ are bounded. Composing with the remaining $C^2$
maps, we obtain
\begin{equation}\label{eq:d/dt(pi(w(t)))}
\Pi'(t)=O(1)
\qquad\text{ and }\qquad
\Pi''(t)=O(1).
\end{equation}

On the other hand, by \eqref{eq:t_w}, \eqref{eq:derivative_t_w} and \eqref{eq:derivative_w(t)}, we have
$
\frac{c_1}{w} \leq \abs{\frac{dt}{dw}(w)} \leq \frac{c_2}{w}
$
and
$
\frac{c_3}{w} \leq \abs{\frac{d^2t}{dw^2}(w)} \leq \frac{c_4}{w}.
$
Similarly, we also deduce
$
C_1w\leq \abs{\frac{dw}{dt}(t)} \leq C_2w
$
and
$
C_1w\leq \abs{\frac{d^2w}{dt^2}(t)} = \abs{-\frac{d^2t/dw^2}{(dt/dw)^3}}  \leq C_2w.
$
Therefore,
\begin{equation}\label{eq:pi'(w)}
\pi'(w)
=
\frac{\Pi'(t)}{\frac{dw}{dt}}
\qquad\text{ and }\qquad
\pi''(w)
=
\frac{
\displaystyle
\Pi''(t)\frac{dw}{dt}
-
\Pi'(t)\frac{d^2w}{dt^2}
}{
\left(\frac{dw}{dt}\right)^3.
}
\end{equation}
Using the estimates above, we obtain
\begin{equation}\label{eq:pi''(w)}
|\pi''(w)|\leq \frac{C}{w^2},
\end{equation}
and combining \eqref{eq:bounds_derivative_pi}, \eqref{eq:d/dt(pi(w(t)))} and \eqref{eq:pi''(w)}, we get
\[
|\psi_{J_k}''(z)|
=
\frac{|\pi''(w)|}{|\pi'(w)|^3}
\leq Cw.
\]
For $w\in I_k^\gamma=[w_k^{\min},w_k^{\max}]$, equation
\eqref{eq:bounds_w_j_k} gives
$
w\leq w_k^{\max}
\leq \frac{C_{\max}}{C_{\min}}w_k^{\min},
$
while \eqref{eq:derivative_psi_z} gives
$
\inf_{z\in J_k}|\psi_{J_k}'(z)|
\geq \frac{w_k^{\min}}{K_6},
$
and we deduce that
$
|\psi_{J_k}''(z)|
\leq
D\inf_{u\in J_k}|\psi_{J_k}'(u)|
$
for some $D>0$ independent of $k$. By the Mean Value Theorem
we then get
\[
\abs{\psi_{J_k}'(x)-\psi_{J_k}'(y)}
\leq
D\cdot \inf_{u\in J_k} \abs{\psi_{J_k}'(u)}\cdot|x-y|,
\]
which is \ref{conf:c6} with $D>0$ and $\alpha=1$.
\end{proof}

Since $\Psi$ is a CIFS, we have the property relating its attractor set and the invariant set of the first-return.

\begin{proposition}[{\cite[Proposition~8]{CunhaNovaesPonce2024HausdorffDimensionCantor}}]\label{ppst:invariant_attractor_equal}
The invariant set $\Lambda_U$ as given in \eqref{eq:local} and the attractor set $\Delta_U$, as given in \eqref{eq:atractor_set}, of the CIFS $\Psi$ are equal.
\end{proposition}
This allows us to use the theory of CIFS to deduce information about the first-return and its invariant set.

\subsection{Proof of Theorem \ref{thm:sliding_bykov}}
We now prove each item of Theorem~\ref{thm:sliding_bykov}. Since $\Lambda = \Delta$, we can use the theory of CIFS to deduce the items of the theorem.

\begin{enumerate}

\item[\ref{thm:symbolic_dynamics}: ] First, we recall that the shift map $\sigma:\mathcal{J}^\N\to \mathcal{J}^\N$ is given by $\sigma(J_1, J_2, J_3, \ldots) := (J_2, J_3, J_4, \ldots)$, and the sequence space $\mathcal{J}^\N$ is endowed with the metric $\opname{dist}:\mathcal{J}^\N\times \mathcal{J}^\N \to [0,+\infty)$ given by
\[
\opname{dist}(a,b) = \opname{dist}((a_n)_{n\ge1}, (b_n)_{n\ge1}) := \begin{cases}
	0, &\quad\text{if } a=b; \\
	1/\kappa, &\quad\text{if } a\neq b,
	\end{cases}
\]
where $\kappa:=\min\set{j\in\N : a(j)\neq b(j)}$. In this metric space $\sigma$ is a continuous function.

We are going to prove that $\proj\circ\,\sigma = \pi\circ\proj$, for $\proj$ defined as in \eqref{eq:proj}. Since $\psi_J = \pi_J^{-1}$, for a given $\omega = (J_1, J_2, J_3, \ldots)\in\cal{J}^{\N}$ we have
\[
\proj \circ\ \sigma(\omega)
= \proj(J_2, J_3, J_4, \ldots) = \psi_{(J_2, J_3, J_4, \ldots)}(\gamma).
\]
On the other hand, we also have
\[
\pi\circ\proj(\omega) = \pi\circ\psi_{(J_1, J_2, J_3, \ldots)}(\gamma) 
= \pi \circ \psi_{J_1}\circ\ \psi_{(J_2, J_3, J_4, \ldots)}(\gamma)
= \psi_{(J_2, J_3, J_4, \ldots)}(\gamma),
\]
since $\psi_{J_1} = \pi_{J_1}^{-1}$. Therefore $\proj\circ\,\sigma = \pi\circ\proj$, and it remains to prove that $\proj$ is a homeomorphism. 

For that, taking arbitrary $a, b\in\cal{J}^{\N}$, with $a\neq b$ and $\opname{dist}(a,b)<\epsilon$, there exists $\kappa=\kappa(\epsilon)\ge1$, such that $a(i) =  b(i)$, for all $1\le i\le\kappa$, and then we have
\[
\proj(a) = \psi_a(\gamma) = \psi_{a(1)}\circ\ldots\circ\psi_{a(\kappa)}\left(\psi_{\sigma^{\kappa+1}(a)}(\gamma)\right)
\]
and
\begin{align*}
\proj(b) = \psi_{ b(1)}\circ\ldots\circ\psi_{ b(\kappa)}\left(\psi_{\sigma^{\kappa+1}(b)}(\gamma)\right)
= \psi_{a(1)}\circ\ldots\circ\psi_{a(\kappa)}\left(\psi_{\sigma^{\kappa+1}(b)}(\gamma)\right),
\end{align*}
and so
\[
	\Bigabs{\proj(a) - \proj(b)} \le \sup_{x,y\in\gamma}\Bigabs{\psi_{(a(1),\ldots,a(\kappa))}(x) - \psi_{(a(1),\ldots,a(\kappa))}(y)}
	\le s^{\kappa}\cdot\sup_{x,y\in\gamma}\abs{x-y}\le s^{\kappa}\cdot\ell(\gamma),
\]
where $\ell(\gamma)$ is the length of $\gamma$, proving that $\proj$ is continuous.

Similarly, for a fixed $x\in\Lambda$, there exists $\epsilon>0$ such that for any $y\neq x$ and $\abs{x-y}<\epsilon$, there exists $\kappa = \kappa(\epsilon)\ge1$ such that $x,y\in\psi_{J_1}\circ\ldots\circ\psi_{J_\kappa}(\gamma)$. Therefore, we have
	$
	\proj^{-1}(x) = (J_1,\ldots,J_{\kappa}, J_{\kappa+1}^x, J_{\kappa+2}^x,\ldots)
	$
	and also 
	$
	\proj^{-1}(y) = (J_1,\ldots,J_{\kappa}, J_{\kappa+1}^y, J_{\kappa+2}^y,\ldots),
	$
	for some indexes $\set{J_i}_{1\le i\le\kappa}, \set{J_i^x,J_i^y}_{i\ge\kappa+1}\subset\cal{J}$.
	From this, we deduce that $\opname{dist}\left(\proj^{-1}(x), \proj^{-1}(y)\right) \le 1/\kappa$, proving that $\proj^{-1}$ is also continuous.
	
	Since trivially $\left(\mathcal{J}^{\N}, \sigma\right) \cong \left(\N^\N, \sigma\right)$, we have the desired result.

	\item[\ref{thm:fractal_geometry}: ] The decomposition can be proved by applying Lemma~\ref{lemma:closure_attractor} with $K_\infty = \set{q}$ to our context, which results in $\ov{\Lambda} = \Lambda \cup \left(\cup_{k\ge 1}\pi^{-k}(q)\right)$. The disjointness comes from the fact that $\Lambda$ is $\pi$-invariant, and $q\notin\Lambda$.

For the bounds on the Hausdorff dimension, we choose $r>0$ such that
$
\sum_{k\geq 1}\frac{C_{\max}}{K_5}
\exp\left(-\frac{2\pi k}{|\kappa|}\right) < 1,
$
and then by using the auxiliary pressure of the system
$
P_1(s) := \sum_{k\geq 1}\left(\frac{C_{\max}}{K_5}
\exp\left(-\frac{2\pi k}{|\kappa|}\right)\right)^s
$
and by the fact that $P$ and $P_1$ are strictly decreasing and continuous on $(0,+\infty)$, with $\lim_{s\to0^+} P(s) = +\infty$ and $P(s)\leq\log P_1(s)$, then $P(1)\leq \log P_1(1) < 0$, so system is regular and the solution for $P(\tilde{s})=0$ satisfies $0<\tilde{s}<1$; by Proposition~\ref{proposition:pressure_hausdorff_dimension} we deduce $0<\hausdorff{\Lambda_U}<1$. Also, this directly implies that $\opname{m}_1(\Lambda_U) = 0$ (see \cite[Theorem~2.1]{Schleicher2007HausdorffDimensionIts}). By applying the identity $\hausdorff{\cup_{k\geq0}X_k} = \sup_{k\geq0}\hausdorff{X_k}$ (see \cite[Theorem~2(2)]{Schleicher2007HausdorffDimensionIts}) in the decomposition just proven, we get $\hausdorff{\ov{\Lambda_U}} = \hausdorff{\Lambda_U}$, and therefore $\opname{m}_1(\ov{\Lambda_U})=0$.

Finally, to prove that $\cl{\Lambda}$ is a Cantor set is to show that it is a compact, totally disconnected set and that all of its points are accumulation points. Since it is the closure of a bounded set, it is compact, and the total disconnectedness comes from the fact that its Hausdorff dimension is less than $1$ (see \cite[Proposition~2.5]{Falconer2006FractalGeometryMathematical}). Finally, take $x = \proj(\omega)$, for some infinite sequence $\omega = (J_1, J_2, J_3,\ldots)\in\mathcal{J}^{\N}$, and consider the sequence $\set{x_j}_{j\geq1}\subset\Lambda$, such that $x_k\in\psi_{(J_1,\ldots,J_k)}(\gamma)$ but $x_k\in\psi_{(J_1,\ldots,J_k, J_{k+1})}(\gamma)$, for $k\geq1$. Then $\set{x_j}_{j\geq1}\subset\Lambda\setminus\set{x}$ and $x_j\to x$; since $x$ is an arbitrary point of $\Lambda$, we have that all points of $\Lambda$ are accumulation points, and since the closure of a set does not add isolated points, we have that $\ov{\Lambda}$ also has this property; therefore $\ov{\Lambda}$ is a Cantor set.

\item[\ref{thm:ergodic_properties}: ] This comes from the fact that $\Psi$ is a regular CIFS (as shown in the proof of Theorem~\ref{thm:fractal_geometry}), and then applying Proposition~\ref{ppst:existence_conformal_measure} and Proposition~\ref{ppst:ergodic_measure}.
\end{enumerate}
\qed
\bigskip

\begin{remark}
For a given $c_s \neq 0$, if $c_u = \frac{M-\alpha}{Sc_s}$, where $\alpha \in \left(-\sqrt{c_s^2+1}, \sqrt{c_s^2+1}\right) \setminus \set{M, M-c_s^2}$, then 
$
c_u\not\in\set{0, c_s/S},
$
satisfying $c_u\neq0$ and condition \ref{cond:non_resonance_undistorted}. Also,
$
\abs{M-Sc_sc_u} = \abs{\alpha} < \sqrt{c_s^2+1} < \sqrt{(c_s^2+1)(c_u^2+1)},
$
satisfying \ref{cond:distortion_bounds}, so with this pair of parameters, the hypotheses of Theorem~\ref{thm:sliding_bykov} are satisfied.
\end{remark}

\begin{remark}\label{rem:sharper_bounds}
In particular, the theory of CIFS provides tools for obtaining sharper estimates for the Hausdorff dimension. By applying Proposition~\ref{ppst:ifs_hausdorff} and Theorem~\ref{thm:cifs_approximation_hausdorff}, for instance, we have
\[
\hausdorff{\Lambda} = \lim_{n\to\infty} \hausdorff{\Lambda_n},
\]
where $\Lambda_n$ is the invariant set of $\pi|_{W^{\gamma}_n}:W^{\gamma}_n\to\gamma$, with $W^{\gamma}_n := \cup_{1\le k\le n}I_k^{\gamma}$. From this we define $\Psi_n := \set{\psi_{J_k}}_{1\le k\le n}$,
$
A:=C_{\min}/K_6,
B:=C_{\max}/K_5,
c:=2\pi/\abs{\kappa},
$
and since
$
\sum_{k=1}^{\infty}e^{-cks}
=
\frac{1}{e^{cs}-1},
$
we obtain
$
\
0 < d_1\leq\hausdorff{\Lambda}\leq d_2,
$
where $d_1$ and $d_2$ are the positive solutions of
\[
\sum_{k=1}^{\infty}A^{d_1}e^{-ckd_1}
= 
\frac{A^{d_1}}{e^{cd_1}-1}=1
\qquad\text{ and }\qquad
\sum_{k=1}^{\infty}B^{d_2}e^{-ckd_2}
=
\frac{B^{d_2}}{e^{cd_2}-1}=1.
\]
\end{remark}

\section{Sliding Bykov Cycles: Examples}\label{section:sliding_bykov_examples}

Now we provide a class of examples of piecewise smooth vector fields, even more piecewise linear vector fields, possessing a sliding Bykov cycle.

\begin{example}[Piecewise Linear Sliding Bykov Cycle]
\label{example:sliding_bykov}

We construct an explicit piecewise linear vector field
$Z=(X,Y)\in\Omega^3(\mathbb{R}^3)$ with switching set
$\Sigma=f^{-1}(0)$, where $f(x,y,z)=z$. Define the smooth vector
fields $X$ in $R^+=\{z>0\}$ and $Y$ in $R^-=\{z<0\}$ by
\begin{equation*}\label{eq:sliding_bykov}
\begin{cases}
X(p)=A_s(p-p_s), & p\in R^+,\\[2mm]
Y(p)=A_Yp+\begin{pmatrix}0\\0\\5\end{pmatrix},
& p\in R^-,
\end{cases}
\end{equation*}
where $p_s=(0,0,1)^T\in R^+$ is the equilibrium of $X$, with
\[
A_s=
\begin{pmatrix}
-1&-2&0\\
2&-1&0\\
0&5&5
\end{pmatrix}
\qquad\text{and}\qquad
A_Y=
\begin{pmatrix}
3&0&0\\
0&3&0\\
0&0&0
\end{pmatrix}.
\]

We now verify that the requirements \ref{SB1}, \ref{SB2}, and
\ref{SB3} are satisfied.

\begin{enumerate}

\item[\ref{SB1}]
First,
$
DX(p_s)=A_s
$
has eigenvalues $\set{5, -1\pm2i}$, so $p_s\in R^+$ is a hyperbolic saddle-focus for $X$.

To determine the sliding region $\Sigma^s\subset\Sigma$, we evaluate
$Xf$ and $Yf$ on $z=0$. Since $f(x,y,z)=z$, we have that
$
Xf(x,y,0)=5(y-1),
$
and
$
Yf(x,y,0)=5,
$
so
$
Xf(x,y,0)<0<Yf(x,y,0)\iff y<1$. Therefore
\[
\Sigma^s
=
\{(x,y,0)^T\in\Sigma:y<1\}
\qquad\text{ and }\qquad
\partial\Sigma^s
=
\{(x,1,0)^T:x\in\mathbb{R}\}.
\]

The Filippov sliding vector field on $\Sigma^s$ is given by
\begin{align*}
Z_\Sigma(x,y)
&=
\left.
\frac{(Yf)X-(Xf)Y}{(Yf)-(Xf)}
\right|_{z=0}
=
\frac{1}{2-y}
\left[
\begin{pmatrix}
-x-2y\\
2x-y
\end{pmatrix}
+
(1-y)
\begin{pmatrix}
3x\\
3y
\end{pmatrix}
\right],
\end{align*}
or equivalently,
\begin{equation*}\label{eq:sliding_field_explicit}
Z_\Sigma(x,y)
=
\frac{1}{2-y}
\begin{pmatrix}
2x-2y-3xy\\
2x+2y-3y^2
\end{pmatrix}.
\end{equation*}

At $p_u=(0,0,0)^T\in\Sigma^s$ we have
$Z_\Sigma(p_u)=0$, so $p_u$ is a pseudo-equilibrium. Besides, 
$
DZ_\Sigma(p_u)
$
has eigenvalues $\set{1\pm i}$, therefore $p_u$ is a hyperbolic pseudo-saddle-focus on $\Sigma^s$.
\item[\ref{SB2}]
The stable eigenspace of $A_s$ corresponds to the eigenvalues
$-1\pm2i$. The invariant stable plane $W^s(p_s)$ is therefore
\[
W^s(p_s)
=
\left\{
(x,y,z):
z=1-\frac14x-\frac34y
\right\}.
\]
Its intersection with
$\partial\Sigma^s=\{z=0,y=1\}$ is the unique point
$
q=(1,1,0)^T.
$

At $q$ we have $(Xf)(q)=0$, while
$
(X^2f)(q)
=
\langle X(q),\nabla(Xf)(q)\rangle
=
5>0,
$
making $q$ a visible fold point.

Now we must prove that $q\in W^u(p_u)$. Since $2-y>0$ for $y\leq 1$, the vector fields $Z_\Sigma$ and
$\tilde{Z}_{\Sigma}:=(2-y)Z_\Sigma$ have the same trajectories with the same orientation in this region, so we instead analyze the system
\[
\begin{cases}
\dot x = 2x-2y-3xy, \\
\dot y = 2x+2y-3y^2.
\end{cases}
\]

Writing this in polar coordinates
$
x=r\cos\theta$
and
$y=r\sin\theta,
$
we obtain
\[
\begin{cases}
\dot r=r(2-3r\sin\theta),\\
\dot\theta=2.
\end{cases}
\]
In particular, we may use $\theta$ as the variable. By setting
$
u=1/r,
$
we get
$
\frac{du}{d\theta}+u=\frac{3}{2}\sin\theta,
$
and since $q=(1,1)$ corresponds to the initial conditions
$
\theta_0=\frac{\pi}{4}$ and $u(\theta_0)=\frac1{\sqrt2},
$
the solution is given by
\[
u(\theta)
=
\frac{e^{\pi/4-\theta}}{\sqrt2}
+\frac34(\sin\theta-\cos\theta).
\]

This implies that $q\in W^u(p_u)$, since
$
u(\theta)\rightarrow\infty$
 as $\theta\rightarrow-\infty
\implies
r(\theta)=1/u(\theta)\rightarrow0$ as $\theta\rightarrow-\infty.
$
Since $\theta=2t+\pi/4$, this means that
$
(x(t),y(t))\rightarrow(0,0)=p_u
\text{ as }t\rightarrow-\infty,
$
therefore $q\in W^u(p_u)$.

It remains to show that the past trajectory of $q$ is contained in
$\{y<1\}$. Since
$
y=r\sin\theta=(\sin\theta)/u,
$
it suffices to prove that $u>0$ for $\theta\leq 0$ and $u>\sin\theta$ for
$0<\theta<\pi/4$.

For $\theta\leq0$, since $\sin\theta - \cos\theta \geq -\sqrt{2}$, we have
$
u(\theta)
\geq
\frac{e^{\pi/4}}{\sqrt2}-\frac{3\sqrt2}{4}
>0.
$

For $0<\theta<\pi/4$, we define
\[
v(\theta):=u(\theta)-\sin\theta
=
\frac{e^{\pi/4-\theta}}{\sqrt2}
-\frac14\sin\theta-\frac34\cos\theta.
\]
We have $v(\pi/4)=0$, while
$
v'(\theta)
=
-\frac{e^{\pi/4-\theta}}{\sqrt2}
-\frac{1}{4}\cos\theta
+\frac{3}{4}\sin\theta.
$
Since $0<\theta<\pi/4$, we have
$
-\frac{1}{4}\cos\theta + \frac{3}{4}\sin\theta
<
\frac{1}{2\sqrt2}
$
and
$
-\frac{e^{\pi/4-\theta}}{\sqrt2}
<
-\frac{1}{\sqrt2},
$
so we can deduce that
$
v'(\theta)<-\frac{1}{2\sqrt2}<0,
$
which implies
$
v(\theta)>v(\pi/4) = 0$, therefore $u(\theta) > \sin(\theta)$
for $0<\theta<\frac{\pi}{4}.
$

Finally, since
$q\in W^s(p_s)$, the orbit continues from the fold into
$
\Gamma_{q\to p_s}\subset R^+\cap W^s(p_s).
$
\item[\ref{SB3}]
The unstable eigenvector of $A_s$ corresponding to the eigenvalue $5$
is given by
$
v_u=(0,0,1)^T,
$
so the unstable manifold of $p_s$ is the $z$-axis through $p_s$. Inside $R^+$ it is parametrized as
\[
\Gamma_{p_s\to p_u}
=
\set{(0,0,z)^T:0<z<1}.
\]
This readily implies
$
\overline{\Gamma_{p_s\to p_u}}\cap\Sigma
= \overline{\Gamma_{p_s\to p_u}}\cap\set{z=0}
=
\{(0,0,0)^T\}
=
\{p_u\},
$
as required.

\end{enumerate}
\end{example}

\begin{remark}
While Example \ref{example:sliding_bykov} does possess a sliding Bykov cycle, we should still verify if it satisfies the hypotheses of Theorem~\ref{thm:sliding_bykov}.

The verification of \ref{cond:linearization} is quite direct, since $X,Z_{\Sigma}\in C^{\infty}$ near their respective equilibria, and $\lambda_s\neq 2a_s$.

For \ref{cond:non_resonance_undistorted}, since $c_s = -1/2$ and $c_u = 1$ we have $c_s - Sc_u\neq0$, regardless of any option of $S\in\set{-1,1}$.

At last, we calculate $L$ in order to show that it satisfies \ref{cond:distortion_bounds}.

We have that $X$ is linear and its local stable manifold is given by
$
W_{loc}^s(p_s)
=
\left\{
z=1-\frac{x}{4}-\frac{3y}{4}
\right\}.
$
Introducing the local coordinates
$
h_s(x,y,z)
=
\left(
x,y,z-1+\frac{x}{4}+\frac{3y}{4}
\right),
$
so that $h_s(p_s)=0$ and $h_s(W^s(p_s))=\{\overline z=0\}$ we have that its inverse is given by
$
h_s^{-1}(\overline x,\overline y,\overline z)
=
\left(
\overline x,
\overline y,
1+\overline z-\frac{\overline x}{4}
-\frac{3\overline y}{4}
\right),
$
therefore
\[
D(h_s^{-1})
=
\begin{pmatrix}
1&0&0\\
0&1&0\\
-\frac14&-\frac34&1
\end{pmatrix}.
\]

We now describe the projection map $\mu_{proj}$.
Along the orbit through
$
u_{out}=(0,0,z_{out}),
$
we have $x(t)=y(t)=0$ and
$
z(t)
=
1+e^{5t}
\left[
z_0-1+
5\int_0^t
e^{-6s}(x_0\sin 2s+y_0\cos 2s)\,ds
\right],
$
so the time of the first intersection with $\Sigma=\{z=0\}$ is
$
t_0=\frac15\log\frac1{1-z_{out}}.
$

The $(x,y)$-component of the flow is
\[
\begin{pmatrix}x(t)\\y(t)\end{pmatrix}
=
e^{-t}
\begin{pmatrix}
\cos 2t&-\sin 2t\\
\sin 2t&\cos 2t
\end{pmatrix}
\begin{pmatrix}x_0\\y_0\end{pmatrix}.
\]
Therefore, linearizing the projection at $u_{out}$ and evaluating at
the corresponding flight time $t_0$ gives
\[
D\mu_{proj}(u_{out})
=
\begin{pmatrix}
e^{-t_0}\cos(2t_0)&-e^{-t_0}\sin(2t_0)&0\\
e^{-t_0}\sin(2t_0)&e^{-t_0}\cos(2t_0)&0\\
0&0&0
\end{pmatrix}.
\]

So we have
\[
D\mu_{proj}(u_{out})
D(h_s^{-1})(\hat u_{out})
\equiv
e^{-t_0}
\begin{pmatrix}
\cos(2t_0)&-\sin(2t_0)\\
\sin(2t_0)&\cos(2t_0)
\end{pmatrix}.
\]

Now we analyze the linearizing map $h_u$. The sliding vector field is
$
Z_\Sigma(x,y)
=
\frac1{2-y}
\begin{pmatrix}
2x-2y-3xy\\
2x+2y-3y^2
\end{pmatrix},
$
and $p_u=(0,0)$ is its pseudo-equilibrium. Again, we have
$
A_u:=DZ_\Sigma(p_u)
=
\begin{pmatrix}
1&-1\\
1&1
\end{pmatrix}.
$
We consider a local diffeomorphism $h_u$ with
$h_u(p_u)=0$ that linearizes $Z_\Sigma$ at $p_u$. In other words, if
$\varphi_{Z_\Sigma}^t$ denotes the flow of the sliding vector field,
then it must satisfy
$
h_u\circ\varphi_{Z_\Sigma}^t
=
e^{tA_u}\circ h_u
$
for all sufficiently small $t$, in all points sufficiently close to
$p_u$.

Since $h_u$ is a local diffeomorphism, $B:=Dh_u(p_u)$ is invertible. By differentiating at $p_u$
the relation
$
h_u\circ\varphi_{Z_\Sigma}^t
=
e^{tA_u}\circ h_u
$
we obtain
$
B e^{tA_u}=e^{tA_u}B.
$
Now, differentiating this identity with respect to $t$ at $t=0$ results in
$
BA_u=A_uB,
$
allowing us to conclude that since $B$ commutes with $A_u$, so does $B^{-1}$.

So we define
$
\widetilde h_u:=B^{-1}\circ h_u,
$
which is also a linearizing map for the same sliding vector field. Furthermore,
$
D\widetilde h_u(p_u)
=
B^{-1}Dh_u(p_u)
=
B^{-1}B
=
\operatorname{Id}_2.
$
Thus we may assume
without loss of generality that
$
Dh_u(p_u)=\operatorname{Id}_2.
$

We then have
$
L
=
e^{-t_0}
\begin{pmatrix}
\cos(2t_0)&-\sin(2t_0)\\
\sin(2t_0)&\cos(2t_0)
\end{pmatrix},
$
implying that both singular values of $L$ are equal to
$
\sigma_1=\sigma_2=e^{-t_0}$, therefore
$
M
=
\frac{\sigma_1^2+\sigma_2^2}{2\sigma_1\sigma_2}
=
1.
$
Also,
$
\det L
=
e^{-2t_0}>0,
$
so
$
S=\operatorname{sign}(\det L)=1.
$

Applying all this in \ref{cond:distortion_bounds}, we have
\[
\abs{M - S c_s c_u} =  \frac{3}{2} < \frac{\sqrt{10}}{2} = \sqrt{(c_s^2 + 1)(c_u^2 + 1)}.
\]

Therefore, Example~\ref{example:sliding_bykov} possesses a sliding Bykov cycle and satisfies the hypotheses of Theorem~\ref{thm:sliding_bykov}, from which we can draw the conclusions.

\end{remark}

\section{Proof of Theorem C}
\label{section:foliations}

We divide the proof of Theorem~\ref{thm:versal_unfolding} into four
steps. The purpose of each step is the following.

\begin{enumerate}
    \item[\textit{Step 1.}]
    \textit{Regularization of the triple crossing.}
    We consider piecewise constant vector fields whose switching set is
    \[
        \Sigma=\{(x,y,z)\in\mathbb{R}^{3}:xyz=0\}.
    \]
    Starting from the eight constant vector fields associated with the
    octants determined by $\Sigma$, we compute the convolution
    regularization and obtain an explicit expression for the resulting
    smooth vector field after a directional blow-up.

    \item[\textit{Step 2.}]
    \textit{Construction of the spatial cusp.}
    We choose a particular collection
    \[
        C=\{C_{\widetilde v}\}_{\widetilde v\in\{-1,1\}^{3}}
    \]
    for which the corresponding regularized vector field $Y_C$ has a
    spatial cusp at the origin.

    \item[\textit{Step 3.}]
    \textit{Construction of a versal unfolding.}
    We perturb the collection $C$ by introducing three parameters and
    show that the corresponding family of regularized vector fields
    provides a versal unfolding of the spatial cusp obtained in Step~2.

    \item[\textit{Step 4.}]
    \textit{Existence of Bykov cycles.}
    We apply the result of Barrientos, Ib\'a\~nez, and Rodr\'iguez
    \cite[Theorem~C]{BarrientosIbanezRodriguez2011HeteroclinicCyclesArising}
    to the versal unfolding constructed in Step~3. This gives a
    bifurcation curve along which the regularized vector field exhibits
    a Bykov cycle, completing the proof of
    Theorem~\ref{thm:versal_unfolding}.
\end{enumerate}

We now carry out these four steps.

\subsection{Triple crossings and piecewise constant fields}

The local codimension of a point of the switching set is the number of
independent switching hypersurfaces passing through it, minus one.  More
precisely, in suitable local coordinates, a point of codimension $0$,
$1$, or $2$ belongs, respectively, to a switching set of the form
\[
    \{x=0\},\qquad \{xy=0\},\qquad \{xyz=0\}.
\]
In particular, the origin is a codimension-$2$ switching point of the
set $\Sigma$ considered here.

For
\[
    \widetilde v=(v_1,v_2,v_3)\in\{-1,1\}^{3},
\]
let
\[
    C=\{C_{\widetilde v}\}_{\widetilde v\in\{-1,1\}^{3}},
    \qquad
    C_{\widetilde v}=(a_{\widetilde v},b_{\widetilde v},
    c_{\widetilde v})
\]
be a collection of constant vector fields.  The corresponding
piecewise constant vector field is
\begin{equation*}
    Z^{C}(x,y,z)=
    \sum_{\widetilde v\in\{-1,1\}^{3}}
    \mathbf{1}_{\{v_1x>0,\,v_2y>0,\,v_3z>0\}}
    C_{\widetilde v}.
\end{equation*}
Thus $Z^{C}$ is constant on each of the eight octants determined by
$\Sigma$.

\subsection{Convolution regularization with a smooth mollifier}\label{subsection:convolution_regularization}

Let $\eta:\mathbb{R}\to[0,\infty)$ be the standard even mollifier
\begin{equation*}
    \eta(s)=\frac{1}{I}
    \begin{cases}
      \displaystyle\exp\!\left(-\frac{1}{1-s^{2}}\right),& |s|<1,\\[2mm]
      0,& |s|\geq1,
    \end{cases}
    \qquad
    I=\int_{-1}^{1}\exp\!\left(-\frac{1}{1-s^{2}}\right)\,ds,
\end{equation*}
and define
\[
    m(u_1,u_2,u_3)=\eta(u_1)\eta(u_2)\eta(u_3).
\]
Then $m\in C_c^{\infty}(\mathbb{R}^{3})$, $m\geq0$, and
$\int_{\mathbb{R}^{3}}m(u)\,du=1$.  For $\varepsilon>0$, put
\[
    m_{\varepsilon}(u)=\frac{1}{\varepsilon^{3}}
    m\!\left(\frac{u}{\varepsilon}\right),
    \qquad
    Z^{C}_{\varepsilon}=m_{\varepsilon}*Z^{C}.
\]
The field $Z^{C}_{\varepsilon}$ is smooth for every $\varepsilon>0$.

Introduce the smooth transition function
\begin{equation*}
    \Phi(s)=2\int_{-\infty}^{s}\eta(r)\,dr-1.
\end{equation*}
Since $\eta$ is even, $\Phi$ is odd.  Moreover,
\begin{gather*}
    \Phi(s)=\pm1\qquad\text{for}\qquad \pm s\geq1,\\
    \kappa:=\Phi'(0)=2\eta(0)>0,
    \qquad \Phi(s)=\kappa s+O(s^{3}).
\end{gather*}

Consider the directional blow-up
\[
    (x,y,z)=\varepsilon(\bar x,\bar y,\bar z).
\]
In the blown-up variables, the differential equation determined by
$Z^{C}_{\varepsilon}$ is, after the constant time rescaling
$d\tau=dt/\varepsilon$,
\begin{equation*}
    \frac{d}{d\tau}(\bar x,\bar y,\bar z)
    =Y_C(\bar x,\bar y,\bar z),
\end{equation*}
where
\begin{equation}\label{eq:YC-smooth-mollifier}
    Y_C(\bar x,\bar y,\bar z)
    =\sum_{\widetilde v\in\{-1,1\}^{3}}
      P_{\widetilde v}(\bar x,\bar y,\bar z)C_{\widetilde v},
    \qquad
    P_{\widetilde v}
    =\frac18\prod_{j=1}^{3}
      \bigl(1+v_j\Phi(\bar x_j)\bigr),
\end{equation}
with $(\bar x_1,\bar x_2,\bar x_3)=(\bar x,\bar y,\bar z)$.
Indeed, after changing variables in the convolution integral, its
$j$-th one-dimensional factor is
\[
    \int_{\mathbb{R}}\eta(r)
    \mathbf{1}_{\{v_j(\bar x_j-r)>0\}}\,dr
    =\frac{1+v_j\Phi(\bar x_j)}{2}.
\]
Notice that the nonsmooth averaging kernel
$\frac18\mathbf{1}_{[-1,1]^{3}}$ would give $\Phi(s)=s$ for
$|s|<1$.  Formula~\eqref{eq:YC-smooth-mollifier} is the corresponding
formula for a genuine $C^{\infty}$ mollifier.

\subsection{A spatial cusp}

We now prove the first assertion of Theorem~\ref{thm:versal_unfolding}.  Recall that a
singularity of a vector field in $\mathbb{R}^{3}$ is a spatial cusp if,
after a smooth change of coordinates and a smooth positive
reparametrization of time, the field can be written as
\[
    \bigl(Y,Z,f(X,Y,Z)\bigr),
    \qquad
    j^{1}f(0)=0,
    \qquad
    \frac{\partial^{2}f}{\partial X^{2}}(0)\neq0.
\]

Consider the collection
\begin{equation}\label{eq:cusp-collection}
    C_{(v_1,v_2,v_3)}=
    \begin{pmatrix}
      v_1-v_3\\
      v_1-v_3-v_1v_2\\
      v_2-v_3
    \end{pmatrix}.
\end{equation}
Set
\[
    u=\Phi(\bar x),\qquad v=\Phi(\bar y),\qquad
    w=\Phi(\bar z).
\]
The character identities on $\{-1,1\}^{3}$ applied to
\eqref{eq:YC-smooth-mollifier} give the exact expression
\begin{equation}\label{eq:YC-uvw}
    Y_C=(u-w,\,u-w-uv,\,v-w).
\end{equation}
We use the local coordinates
\begin{equation}\label{eq:cusp-coordinates}
    (X,Y,Z)=(u,\,u-w,\,u-v).
\end{equation}
They are well defined near the origin because $\Phi'(0)=\kappa>0$.
If
\[
    A(s)=\Phi'\bigl(\Phi^{-1}(s)\bigr),
\]
then $A(s)=\kappa+O(s^{2})$, since $\Phi$ is odd.  Transforming the
vector field~\eqref{eq:YC-uvw} by
\eqref{eq:cusp-coordinates}, we obtain
\[
    \dot X=\kappa Y+O(3),\qquad
    \dot Y=\kappa Z+O(3),\qquad
    \dot Z=\kappa X(X-Z)+O(3).
\]
After rescaling time by the positive constant $\kappa$, the $2$-jet is
\begin{equation}\label{eq:nilpotent-two-jet}
    X'=Y+O(3),\qquad
    Y'=Z+O(3),\qquad
    Z'=X(X-Z)+O(3).
\end{equation}
The standard preparation theorem for nilpotent singularities allows us
to choose local coordinates, without changing this $2$-jet, in which
the first two components are exactly $Y$ and $Z$.  The last component
$f$ then satisfies
\[
    j^{1}f(0)=0,
    \qquad
    \frac{\partial^{2}f}{\partial X^{2}}(0)=2\neq0.
\]
Consequently, the origin is a spatial cusp, or equivalently a
three-dimensional nilpotent singularity of codimension three.

\subsection{Versal unfolding and Bykov cycles}

We perturb the collection~\eqref{eq:cusp-collection} using three
parameters.  Namely, let
\begin{equation}\label{eq:perturbed-collection}
\begin{split}
    C_{(v_1,v_2,v_3)}^{a,b,c}
    ={}&C_{(v_1,v_2,v_3)}
       -a\begin{pmatrix}0\\1\\0\end{pmatrix}
       -b\begin{pmatrix}0\\v_1-v_3\\0\end{pmatrix}\\
       &{}-c\begin{pmatrix}0\\v_1-v_2\\0\end{pmatrix},
       \qquad (a,b,c)\in\mathbb{R}^{3}.
\end{split}
\end{equation}
Let $Y_{a,b,c}$ be the vector field obtained from this collection by
\eqref{eq:YC-smooth-mollifier}.  Before applying
\eqref{eq:cusp-coordinates}, the character identities give
\begin{equation*}
    Y_{a,b,c}=
    \bigl(u-w,\,
    u-w-uv-a-b(u-w)-c(u-v),\,v-w\bigr).
\end{equation*}
After the coordinate preparation used above, the last component has
the form
\begin{equation}\label{eq:versal-prepared-form}
\begin{split}
    f_{a,b,c}(X,Y,Z)
    ={}&X(X-Z)+\alpha(a,b,c)+\beta(a,b,c)Y
       +\gamma(a,b,c)Z\\
       &+O\bigl(\abs{(X,Y,Z)}^{3}\bigr),
\end{split}
\end{equation}
where
\begin{equation}\label{eq:versal-transversality}
    \det D(\alpha,\beta,\gamma)(0,0,0)\neq0.
\end{equation}
Indeed, at the level of jets, the three perturbations in
\eqref{eq:perturbed-collection} independently generate the classes of
$1$, $Y$, and $Z$ in the quotient space associated with the nilpotent
singularity.  These classes form a basis of that three-dimensional
quotient.  Condition~\eqref{eq:versal-transversality} is therefore the
transversality condition for a miniversal unfolding.  Hence
$Y_{a,b,c}$ is a versal, and in particular generic, unfolding of the
spatial cusp.

It remains to prove the last assertion of Theorem~\ref{thm:versal_unfolding}.  We use the
following result of Barrientos, Ib\'a\~nez, and Rodr\'iguez
\cite[Theorem~C]{BarrientosIbanezRodriguez2011HeteroclinicCyclesArising}: every generic unfolding
of a three-dimensional nilpotent singularity of codimension three
contains a one-sided bifurcation curve of topological Bykov cycles.
Each such cycle consists of two saddle-focus equilibria with different
stability indices and two heteroclinic connections between them.

Equations~\eqref{eq:nilpotent-two-jet} and
\eqref{eq:versal-transversality} verify precisely the two hypotheses
needed to apply that result: the origin is a nilpotent singularity of
codimension three and the family $Y_{a,b,c}$ is a generic unfolding.
The higher-order terms introduced by the smooth mollifier cause no
additional obstruction.  Under the weighted nilpotent rescaling used
in \cite{BarrientosIbanezRodriguez2011HeteroclinicCyclesArising}, the family
\eqref{eq:versal-prepared-form} takes the form
\[
    \dot x=y,\qquad
    \dot y=z,\qquad
    \dot z=\nu_1+\nu_2y+\nu_3z+x^{2}+O(\delta),
\]
where $\delta>0$ is the rescaling parameter.  Thus the terms depending
on the particular mollifier belong to the $O(\delta)$ remainder covered
by the cited theorem.  It follows that there are $\delta_0>0$ and a
curve
\[
    (a,b,c)=\bigl(a(\delta),b(\delta),c(\delta)\bigr),
    \qquad 0<\delta<\delta_0,
\]
converging to $(0,0,0)$ such that $Y_{a,b,c}$ possesses a topological
Bykov cycle.  In particular, Bykov cycles occur for arbitrarily small
perturbations of the collection~\eqref{eq:cusp-collection}.  This
completes the proof of Theorem~\ref{thm:versal_unfolding}.
\qed

\bigskip

\begin{remark}
No reversibility assumption on the full regularized family is needed.
Reversibility is used in \cite{BarrientosIbanezRodriguez2011HeteroclinicCyclesArising} to locate
the cycle in the limit family; the splitting equations and the implicit
function theorem then continue the cycle to a generic unfolding with
higher-order terms.
\end{remark}

\section{Conclusions}
\label{section:conclusions}
We have described three mechanisms relating Bykov cycles to piecewise smooth dynamics.
First, crossing Bykov cycles resemble the main mechanism of the first-return map from classical smooth
Bykov cycles, with the switching manifold acting only as an interface, and the corresponding
return dynamics contains horseshoes accumulating on the cycle. Second, sliding Bykov
cycles exhibit a genuinely non-smooth mechanism. Under the stated conditions, the associated return map generates a conformal iterated function system, from which we can deduce dynamical and geometric information.

Lastly, a different kind of result arises from the study of regularization. For a suitable piecewise constant vector field with a triple switching set, a convolution regularization produces a smooth
spatial cusp with its versal unfolding containing a classical Bykov cycle.

These results indicate that Bykov-type dynamics interacts with non-smooth systems in two complementary ways: it may be intrinsic to the Filippov dynamics, through crossing or sliding connections, or emerge in the smooth families obtained by regularization.

The sliding case also suggests further questions. In particular, it would be natural to
study the persistence and bifurcations of these cycles under perturbations, as well
as the relation between their invariant sets and the dynamics of regularized
systems near the switching manifold.

\appendix


\section{Asymptotics of Remainder Functions}\label{appendix:asymptotic_of_remainder}
For the sake of completeness, in this section we prove a lemma used in many parts of the main text, allowing us to estimate some error terms of functions and their derivatives.
\begin{lemma}[Taylor remainder along a $C^2$ curve]
\label{lem:taylor_remainder_curve}
Let $G:\mathbb{R}^n\to\mathbb{R}^m$ be of class $C^2$ in a
neighborhood of $u_*$, and let $\beta:[T,\infty)\to\mathbb{R}^n$ be a
$C^2$ curve such that $u_*+\beta(t)$ remains in a sufficiently small neighborhood of $u_*$. Define
\[
R_G(\beta) := G(u_*+\beta)-G(u_*)-DG(u_*)\beta.
\]
Then there exists a constant $C_R>0$ such that
\[
\norm{R_G(\beta(t))} \le C_R\norm{\beta(t)}^2,
\]
\[
\norm{\frac{d}{dt}R_G(\beta(t))}
= \norm{\bigl(DG(u_*+\beta(t))-DG(u_*)\bigr)\beta'(t)}
\le C_R\norm{\beta(t)}\norm{\beta'(t)},
\]
and
\[
\norm{\frac{d^2}{dt^2}R_G(\beta(t))}
\le C_R\left(
\norm{\beta'(t)}^2+
\norm{\beta(t)}\norm{\beta''(t)}
\right).
\]
\end{lemma}
\begin{proof}
The first estimate is the usual Taylor estimate for $C^2$ functions.
For the derivative, differentiating $R_G$ along the
curve $\beta(t)$ results, by the chain rule, in
\[
\frac{d}{dt}R_G(\beta(t))
= DG(u_*+\beta(t))\beta'(t)-DG(u_*)\beta'(t)
= \bigl(DG(u_*+\beta(t))-DG(u_*)\bigr)\beta'(t).
\]
Since $DG$ is locally Lipschitz, we have
\[
\norm{DG(u_*+\beta(t))-DG(u_*)} \le C_R\norm{\beta(t)},
\]
which results in
\[
\norm{\frac{d}{dt}R_G(\beta(t))}
\le C_R\norm{\beta(t)}\norm{\beta'(t)}.
\]
For the second derivative, differentiating once more and using the chain
rule gives
\[
\norm{\frac{d^2}{dt^2}R_G(\beta(t))}
\leq
C_R\norm{\beta'(t)}^2\\
+
\norm{DG(u_*+\beta(t))-DG(u_*)}
\norm{\beta''(t)}.
\]
Since $D^2G$ is bounded and $DG$ is locally Lipschitz,
\[
\norm{\frac{d^2}{dt^2}R_G(\beta(t))}
\leq C_R\left(
\norm{\beta'(t)}^2+
\norm{\beta(t)}\norm{\beta''(t)}
\right).
\]
\end{proof}
In practice, if $\beta(t)=O(g(t))$,
$\beta'(t)=O(g'(t))$, and $\beta''(t)=O(g''(t))$, then the above
estimates give us
\[
R_G(\beta(t))=O(g(t)^2),
\quad
\frac{d}{dt}R_G(\beta(t))=O(g(t)g'(t)),
\quad\text{ and }\quad
\frac{d^2}{dt^2}R_G(\beta(t))
=O\left(g'(t)^2+g(t)g''(t)\right).
\]

\section{Transversality of $\mu$ and $Z_{\Sigma}$}\label{appendix:transversality}
We now prove that the ($C^2$-linearized) projected curve $\hat{\mu}(t)$ is transversal to the ($C^2$-linearized) sliding vector field $\hat{Z}_{\Sigma}$ in a ($C^2$-linearized) planar neighborhood of $p_u$.

By differentiating \eqref{eq:projected_curve_with_error} with respect to $t$, and using Lemma~\ref{lem:taylor_remainder_curve} to control the derivative of the Taylor remainder, we obtain
\begin{equation}\label{eq:projected_derivative_linear}
\hat{\mu}'(t)
=
r_{in}e^{-a_st}
\left[
-a_sL\hat{v}_0(t)
+b_sLJ\hat{v}_0(t)
+O\!\left(e^{-\delta t}\right)
\right],
\end{equation}
where $\delta=\min\set{a_s,\lambda_s}$.

On the other hand, the linearized sliding vector field $\hat{Z}_{\Sigma}$ evaluated on $\hat{\mu}(t)$ is given by
\begin{equation}\label{eq:sliding_field_linear}
\hat{Z}_{\Sigma}(\hat{\mu}(t)) = \hat{A}_u \hat{\mu}(t) = r_{in} e^{-a_s t} [a_u L \hat{v}_0(t) + b_u J L \hat{v}_0(t) + O(e^{-\delta t})],
\end{equation}
where $\hat{A}_u := a_u I + b_u J$ and $J := R_{\pi/2} =
	\begin{pmatrix}
	0 & -1 \\ 1 & 0
	\end{pmatrix}$.

The tangencies between $\hat{\mu}(t)$ and $\hat{Z}_{\Sigma}$ are then determined by the roots of the determinant
\begin{equation}\label{eq:transversality_determinant}
F(t) =\det[\hat{\mu}'(t), \hat{Z}_{\Sigma}(\hat{\mu}(t))].
\end{equation}
Substituting~\eqref{eq:projected_derivative_linear} and~\eqref{eq:sliding_field_linear} into the determinant gives
\begin{align*}\label{eq:transversality_determinant}
F(t) = r_{in}^2 e^{-2a_s t} \cdot \det\left[-a_s L\hat{v}_0(t) + b_s LJ\hat{v}_0(t) + O(e^{-\delta t}), \, a_u L\hat{v}_0(t) + b_u JL\hat{v}_0(t) + O(e^{-\delta t})\right]
\end{align*}
which can be expanded as
\begin{align*}
r_{in}^2 e^{-2 a_s t} \big[ &\det(-a_s L \hat{v}_0(t), a_u L \hat{v}_0(t)) + \det(-a_s L \hat{v}_0(t), b_u J L \hat{v}_0(t)) + \\
&\det(b_s L J \hat{v}_0(t), a_u L \hat{v}_0(t)) + \det(b_s L J \hat{v}_0(t), b_u J L \hat{v}_0(t)) + O(e^{-\delta t}) \big].
\end{align*}

Applying the identities $\det(v, v) = 0$, $\det(v, Jv) = \norm{v}^2$, and $\det(L J v, L v) = -\det L\norm{v}^2$, we can deduce that
\[
\begin{cases}
\det(-a_s L \hat{v}_0(t), a_u L \hat{v}_0(t)) = 0, \\
\det(-a_s L \hat{v}_0(t), b_u J L \hat{v}_0(t)) = -a_s b_u \norm{L \hat{v}_0(t)}^2, \\
\det(b_s L J \hat{v}_0(t), a_u L \hat{v}_0(t)) = -a_u b_s \det L, \\
\det(b_s L J \hat{v}_0(t), b_u J L \hat{v}_0(t)) = b_s b_u \det(L J \hat{v}_0(t), J L \hat{v}_0(t)),
\end{cases}
\]
which gives us
\begin{equation}\label{eq:F_reduced}
F(t) = r_{in}^2 e^{-2 a_s t} \left[ -a_u b_s \det L - a_s b_u \norm{L \hat{v}_0(t)}^2 + b_s b_u \det(L J \hat{v}_0(t), J L \hat{v}_0(t)) + O(e^{-\delta t}) \right].
\end{equation}

Now we will need the following lemma, where we recall that $\sigma_1\ge\sigma_2 > 0$ are the singular values of $L$, $S := \sgn(\det L)$ and $M := \frac{\sigma_1^2 + \sigma_2^2}{2 \sigma_1 \sigma_2} \ge 1$, as given in \eqref{eq:L_quantities}.

\begin{lemma}\label{lem:wave}
We have
\begin{equation*}
F(t) = r_{in}^2 e^{-2 a_s t} \left[ C + A \cos(2 b_s t - \phi) + O(e^{-\delta t}) \right],
\end{equation*}
where $\phi$ is a phase shift and the coefficients are
\[
C = -b_u b_s \sigma_1 \sigma_2 (M c_s + S c_u) \qquad \text{and} \qquad A = -b_s b_u \sigma_1 \sigma_2 \sqrt{M^2 - 1} \sqrt{c_s^2 + 1}.
\]
\end{lemma}

\begin{proof}[Proof of Lemma~\ref{lem:wave}]
From the expression of $F(t)$ as given in \eqref{eq:F_reduced}, we consider the singular value decomposition
$L=UDV^T$. Here $U,V\in O(2)$ are orthogonal matrices, and $D := \operatorname{diag}(\sigma_1, \sigma_2)$, where $\sigma_1 \ge \sigma_2 > 0$. In particular,
$\det L=\det(U)\det(V)\sigma_1\sigma_2 =S\sigma_1\sigma_2$.

Since $V^T$ is orthogonal and $\norm{\hat{v}_0(t)}=1$, there exists an
angle $\alpha(t)$ such that
\[
\hat{v}_0^V(t) := V^T\hat{v}_0(t)
=
\begin{pmatrix}
\cos\alpha(t)\\
\sin\alpha(t)
\end{pmatrix},
\]
and because $\hat{v}_0(t)$ rotates with angular velocity $b_s$ and
$V^T$ is a fixed orthogonal transformation, we have
$\alpha(t)=\pm b_st+\alpha_0$
for some constant $\alpha_0$.

We then calculate
\begin{equation*}
\norm{L\hat{v}_0(t)}^2
=
\norm{UDV^T\hat{v}_0(t)}^2
=
\norm{D\hat{v}_0^V(t)}^2\\
=
\sigma_1^2\cos^2\alpha(t)
+\sigma_2^2\sin^2\alpha(t)\\
=
\frac{\sigma_1^2+\sigma_2^2}{2}
+
\frac{\sigma_1^2-\sigma_2^2}{2}
\cos(2\alpha(t)).
\end{equation*}

For the determinant in \eqref{eq:F_reduced}, by using the identities
\[
U^TJU=(\det U)J
\qquad\text{ and }\qquad
V^TJV=(\det V)J,
\]
we deduce
\begin{equation*}
\det(LJ\hat{v}_0(t),JL\hat{v}_0(t))
=
-\det(V)
\frac{\sigma_1^2-\sigma_2^2}{2}
\sin(2\alpha(t)).
\end{equation*}

Since
\[
\frac{\sigma_1^2-\sigma_2^2}{2\sigma_1\sigma_2}
=
\sqrt{M^2-1},
\]
we can deduce that the leading term from the determinant of \eqref{eq:F_reduced} can be written as
\begin{equation*}%
-b_sb_u\sigma_1\sigma_2
\Big[
Mc_s+Sc_u
+\sqrt{M^2-1}
\big(
c_s\cos(2\alpha(t))
+\det(V)\sin(2\alpha(t))
\big)
\Big].
\end{equation*}

The trigonometric expression can be written as
\[
B_1\cos(2b_st)+B_2\sin(2b_st),
\]
where $B_1,B_2\in\mathbb{R}$
. By using the identity $B_1\cos(x) + B_2\sin(x) = B\cos(x+\varphi)$, for some
$B\geq0$ and $\varphi\in\mathbb{R}$, we deduce
\[
B_1\cos(2b_st)+B_2\sin(2b_st)
=
B\cos(2b_st-\phi).
\]
Moreover, the amplitude of both trigonometric expressions must match, so we can deduce that $B = \sqrt{B_1^2 + B_2^2} = \sqrt{c_s^2 + \det(V)^2} = \sqrt{c_s^2 + 1} \ge0$, so if we define
\[
A
:=
-b_sb_u\sigma_1\sigma_2
\sqrt{M^2-1}\sqrt{c_s^2+1}
\qquad\text{ and }\qquad
C := -b_sb_u\sigma_1\sigma_2(Mc_s+Sc_u),
\]
then \eqref{eq:F_reduced} becomes
\begin{equation*}%
F(t)
=
r_{in}^2e^{-2a_st}
\left[
C+A\cos(2b_st-\phi_0)
+
O(e^{-\delta t})
\right],
\end{equation*}
for some $\phi_0\in\R$.
\end{proof}

Now, we are in position to prove the main result.

\begin{proposition}\label{pps:F_greater_zero}
For $F(t)$ as in \eqref{eq:transversality_determinant}, there exist $c_F,C_F>0$ and $T>0$ such that
$
0 < c_F e^{-2a_st} \leq \abs{F(t)} \leq C_F e^{-2a_st},
$
for all $t\geq T$. In particular, $\hat{\mu}'(t)$ and $\hat{Z}_{\Sigma}(\hat{\mu}(t))$ are transversal.
\end{proposition}
\begin{proof}[Proof of Proposition~\ref{pps:F_greater_zero}]
In order $\hat{\mu}(t)$ and $\hat{Z}_{\Sigma}$ to be transversal in a punctured neighborhood of $p_u$, we must have $F(t)\neq0$ in it, so we write $F(t)$ in the form given by Lemma~\ref{lem:wave} and impose
\begin{equation}\label{eq:C-A}
\abs{C}-\abs{A} > 0 \iff \abs{M c_s + S c_u} > \sqrt{M^2 - 1} \sqrt{c_s^2 + 1}.
\end{equation}

We can deduce this by squaring both sides of \ref{cond:distortion_bounds}, where we get
\[
M^2 - 2MSc_sc_u+c_s^2c_u^2 < c_s^2c_u^2 + c_s^2 + c_u^2 +1 \implies 0 < c_s^2 + c_u^2 + 1 - M^2+2MSc_sc_u = (Mc_s+Sc_u)^2 - (M^2-1)(c_s^2+1),
\]
where the last expression is equivalent to \eqref{eq:C-A}, as desired.

Therefore, whenever 
\ref{cond:distortion_bounds} is satisfied, the unperturbed wave $C + A \cos(2 b_s t - \phi)$ is bounded away from $0$, and since the error term $O(e^{-\delta t}) \to 0$ as $t \to \infty$, then for all sufficiently large $t$ we have $F(t) \neq 0$. Therefore, we must have $
0 < c_F e^{-2a_st} \leq \abs{F(t)} \leq C_F e^{-2a_st},
$. Also, this implies that $\hat{\mu}(t)$ and $\hat{Z}_{\Sigma}$ are transversal in a punctured neighborhood of $p_u$, which is equivalent to the transversality of $\mu(t)$ and $Z_{\Sigma}$ on a punctured neighborhood of $p_u$ as well.
\end{proof}

\section{Overview of Conformal Iterated Function Systems (CIFS)}\label{appendix:cifs}

This section is devoted to presenting the main tool for the proof of Theorem~\ref{thm:sliding_bykov}. 

Let us consider a compact set $K = \ov{\opname{Int}(K)} \subset \R^d$ and a family of maps $\mathcal{F} := \set{f_i : K\to K}_{i\in I}$, where $I$ is an at most countable index set and such that, for all $x, y \in K$ and for each $i \in I$, there exists a constant $s_i < 1$ such that
\[
\abs{f_i(x) - f_i(y)}\le s_i\cdot\abs{x-y}.
\]
The family $\mathcal{F}$ is called an \emph{iterated function system (IFS)}.

We have a special class of IFS called \emph{conformal IFS (CIFS)}, which are IFS satisfying the following conditions.

\begin{enumerate}[label = (C\arabic*), ref = (C\arabic*)]
	\item\label{conf:c1} For each $i\in I$, the map $f_i$ is an injection of $K$ into itself;
	\item\label{conf:c2} For some $s<1$ it satisfies
	\[
	\abs{f_i(x) - f_i(y)} \leq s\cdot\abs{x-y},
	\]
	for all $i\in I$ and for all $x,y\in K$;
	\item\label{conf:c3} The set $K$ is connected and each $f_i$ satisfies $f_i(\opname{Int}(K))\subset\opname{Int}(K)$ and $f_i(\opname{Int(K)})\cap f_j(\opname{Int(K)}) = \emptyset$, for $i\neq j$;
	\item\label{conf:c4} There is an open set $V\supset K$ such that each $f_i$ can be extended to a $C^{1,\varepsilon}$-diffeomorphism $f_i^V: V\to V$; besides, these extensions are conformal maps, that is, $Df_i^V(x) = \rho_{x,i}\cdot\opname{Isom}_{x,i}$, where $\rho_{x,i}\in\R$ and $\opname{Isom}_{x,i}:\R^d\to\R^d$ is an isometry, for any $x\in V$ and for all $i\in I$;
	\item\label{conf:c5} The following inequality is satisfied.
	\[
	\inf_{x\in\partial K}\inf_{0< \varepsilon< 1}\frac{\opname{m}_d(B_\varepsilon(x)\cap\opname{Int(K)})}{\opname{m}_d(B_\varepsilon(x))} > 0,
	\]
	where $\opname{m}_d$ is the $d$-dimensional Lebesgue measure;
	\item\label{conf:c6} There are $D\geq 1$ and $\alpha>0$ such that, for any $i\in I$, we have
	\[
	\abs{\norm{Df_i^V(x)}_{\operatorname{op}} - \norm{Df_i^V(y)}_{\operatorname{op}}} \leq D\cdot\norm{(Df_i^V)^{-1}}^{-1}_{\infty}\cdot\abs{x-y}^{\alpha},
	\]
	where $f_i^V$ is the map described in \ref{conf:c4}, $x,y\in V$, the norm $\norm{\cdot}_{\operatorname{op}}$ is the operator norm and $\norm{f}_{\infty} := \sup_{z\in V}\set{\norm{f(z)}_{\operatorname{op}}}$.
\end{enumerate}

For $k \ge 1$, the (finite) sequences given by $\eta = (\eta_1,\ldots,\eta_k) \in \bigcup_{j\ge1} I^j$ can be associated to the map
$
f_{\eta} := f_{\eta_1}\circ\cdots\circ f_{\eta_k}.
$
Also, for $\omega \in \cal{J}^k$ and $\eta \in \cal{J}^l$, we define $\omega\eta := (\omega_1,\ldots,\omega_k,\eta_1,\ldots,\eta_l)\in\cal{J}^{k+l}$, inducing the composition $\psi_\omega \circ\psi_\eta$.
For CIFS, we can define the projection map $\proj : I^\N \to K$ where, for every (infinite) sequence $\omega = (\omega_1, \omega_2, \omega_3, \ldots) \in I^\N$, it is given by the intersection of non-empty nested compact sets with diameter approaching $0$, given by
\begin{equation}\label{eq:proj}
\set{\proj(\omega)} = \bigcap_{k\ge1}f_{\omega|_k}(K),
\end{equation}
where $\omega|_k := (\omega_1, \ldots, \omega_k)$ for each $k \ge1$. Notice that in these conditions, $\proj$ is well-defined.

Using these notions, the \emph{attractor set} of an IFS is the set
\begin{equation}\label{eq:atractor_set}
\Delta := \proj\left(I^{\N}\right) = \bigcup_{\omega\in I^{\N}}\set{\proj(\omega)} = \bigcup_{\omega\in I^{\N}} \bigcap_{k\ge1}f_{\omega|_k}(K).
\end{equation}

Now here are some results of the theory of CIFS.

\begin{proposition}[{\cite[Propositions~9.6 and~9.7]{Falconer2006FractalGeometryMathematical}}]\label{ppst:ifs_hausdorff}
For a finite IFS $\set{f_i:K\to K}_{1\leq i\leq k}$ such that, for some $0<b_i\leq c_i < 1$, it satisfies $b_i\abs{x-y}\leq \abs{f_i(x)-f_i(y)}\leq c_i\abs{x-y}$, for all $x,y\in K$, and such that $f_i(\Delta)\cap f_j(\Delta)=\emptyset$. Then $s\leq \hausdorff{\Delta}\leq t$, where $s$ and $t$ are the unique solutions to $\sum_{1\leq i\leq k}b_i^s = 1$ and $\sum_{1\leq i\leq k}c_i^s = 1$.
\end{proposition}

If the IFS is conformal, we have an approximation property.

\begin{theorem}[{\cite[Theorem~3.15]{MauldinUrbanski1996DimensionsMeasuresInfinite}}]\label{thm:cifs_approximation_hausdorff}
For a CIFS $\set{f_i:K\to K}_{i\in I}$ with attractor set $\Delta$, we have
\[
\hausdorff{\Delta} = \sup_{J\in \operatorname{Fin}(I)}\set{\hausdorff{\Delta_J} : \Delta_J \text{ is the attractor set of } \set{f_j}_{j\in J}},
\]
where $\operatorname{Fin}(I)$ is the collection of all finite subsets of $I$.
\end{theorem}

We also have the following characterization of the closure of the attractor set. Consider $I^* := \cup_{n\in\N}I^n$, and $K_{\infty}$ the set of limit points of all sequences $x_i\in f_i(K)$, with $i\in I'$, where $I'$ ranges over all infinite subsets of $I$.

\begin{lemma}[{\cite[Lemma~2.1]{MauldinUrbanski1996DimensionsMeasuresInfinite}}]\label{lemma:closure_attractor}
If $\lim_{i\in I}\opname{diam}(f_i(K_\infty)) = 0$, then $\ov{\Delta} = \Delta \cup \left(\bigcup_{\omega\in I^*}f_{\omega}(K_\infty)\right)$.
\end{lemma}

Now we need to introduce an important concept.
 For each $k\ge1$, consider the auxiliary functions $P_k:[0,+\infty)\to\R\cup\set{+\infty}$ given by
$
P_k(t) := \sum_{\eta\in I^k}\norm{Df_{\eta}}_{\infty}^t.
$ The \emph{pressure of the system} is the function $P:[0,+\infty)\to\R\cup\set{+\infty}$ defined as
$
P(t) := \lim_{k\to+\infty}\frac{1}{k}\log P_k(t)= \lim_{k\to+\infty}\frac{1}{k}\log\left(\sum_{\eta\in I^k}\norm{Df_{\eta}}_{\infty}^t\right).
$

We say that a CIFS is \emph{regular} if its pressure function has a unique zero,  $\hat{t}\ge0$ with $P(\hat{t})=0$, we call the CIFS $\hat{t}$-regular. For regular systems, we have the following result.

\begin{proposition}[{\cite[Theorem 7.4]{Mauldin1995InfiniteIteratedFunction}}]
\label{proposition:pressure_hausdorff_dimension}
Let $\Delta$ be the invariant set of a  $\hat{t}$-regular CIFS. Then, $\hausdorff{\Delta} = \hat{t}$.
\end{proposition}

Finally, we introduce a definition.

\begin{definition}[$t$-Conformal Measures]
\label{definition:t_conformal_measures}Let $t\ge0$ be a real number, then a measure $m_t$ is called a \emph{$t$-conformal measure (for the CIFS $\cal{F}$)} if
\begin{enumerate}
    \item $m_t(\Delta) = 1$;
    \item $m_t(f_i(B)) = \int_{B}\norm{Df_i}^tdm_t$, for all $i\in I$ and for all Borel sets $B\subset\Delta$;
    \item $m_t(f_i(K)\cap f_j(K)) = 0$, for all $i\in I$, $i\neq j$.
\end{enumerate}
\end{definition}

This definition allows us to introduce the two last results of measure theory linked to CIFS.

\begin{proposition}[{\cite[Theorem 7.2]{Mauldin1995InfiniteIteratedFunction}}]
\label{ppst:existence_conformal_measure}
A CIFS $\cal{F}$ is $\hat{t}$-regular if and only if there exists a $\hat{t}$-conformal measure $m_{\hat{t}}$ for $\cal{F}$.
\end{proposition}

We now introduce some definitions and notations: for each $\eta = (\eta_1,\ldots,\eta_k)\in \bigcup_{j} I^j$, we define $\abs{\eta} := k$ as the \emph{length of $\eta$}, and the set $[\eta] := \set{\omega\in I^\N : \omega|_{\abs{\eta}} = \eta}$ as the \emph{cylinder generated by $\eta$}. We then have the following.
\begin{proposition}[{\cite[Lemma 3.6 and Theorem 3.8]{MauldinUrbanski1996DimensionsMeasuresInfinite}}]\label{ppst:measure_on_index_set}
Let $\cal{F}$ be a $\hat{t}$-regular CIFS, and let $m_{\hat{t}}$ be the $\hat{t}$-conformal measure  given by Proposition~\ref{ppst:existence_conformal_measure}.

Then there exists a unique Borel probability measure $\nu$ on $I^\N$ with
$
\nu([\eta]) = \int_K\norm{Df_{\eta}}^{\hat{t}}dm_{\hat{t}},
$
for all $\eta\in\cup_{j\ge 1}I^j$.
Besides, there is a unique measure $\nu^*$ on $I^\N$, with the property of being invariant under the shift $\sigma$ and being equivalent to $\nu$.
\end{proposition}

Considering the measure $\nu^*$ on the coding space associated to $\nu$, and the map $T:\Delta\to\Delta$ defined by $T(\proj(\omega)) := \proj(\sigma(\omega))$, we have the last proposition.

\begin{proposition}[{\cite[Theorem 8.2]{Mauldin1995InfiniteIteratedFunction}}]
\label{ppst:ergodic_measure}
Let $\cal{F}$ be a $\hat{t}$-regular CIFS. Then $m_{\hat{t}} = \nu\circ\proj^{-1}$, and the measure $m_{\hat{t}}^* := \nu^*\circ \proj^{-1}$ on $\Delta$ is the unique invariant and ergodic measure with respect to $T$ and satisfying $m_{\hat{t}}^* \sim m_{\hat{t}}$ on $\Delta$.
\end{proposition}

The construction of all of these measures can be found in detail in \cite[Section 3]{MauldinUrbanski1996DimensionsMeasuresInfinite} and \cite[Section 8]{Mauldin1995InfiniteIteratedFunction}.


\bibliographystyle{plain}
\bibliography{references.bib}

\section*{Acknowledgements}
MGCC is partially supported by S\~{a}o Paulo Research Foundation (FAPESP) grant 2025/01340-7. PRS is partially supported by ANR-23-CE40-0028, 
FAPESP grant 2023/02959-5, CAPES  MATH-Amsud Project
grant: 88881.179491/2025-01,  FAPESP grant 2024/15612-6 and 
CNPq grant 302154/2022-1.

\end{document}